\documentclass[a4paper]{amsart}

\usepackage{amsfonts}
\usepackage{amssymb}
\usepackage{amsthm}
\usepackage{amsmath}
\usepackage[all]{xy}

\SelectTips{cm}{}

\theoremstyle{plain}
\newtheorem{thm}{Theorem}[section]
\newtheorem{prop}[thm]{Proposition}
\newtheorem{cor}[thm]{Corollary}
\newtheorem{lem}[thm]{Lemma}

\theoremstyle{definition}
\newtheorem{defn}[thm]{Definition}

\theoremstyle{remark}
\newtheorem{rem}[thm]{Remark}

\numberwithin{equation}{section}

\newcommand{\sN}{\mathcal{N}}
\newcommand{\sS}{\mathcal{S}}

\newcommand{\bt}{\mathbf{t}}
\newcommand{\bs}{\mathbf{s}}
\newcommand{\br}{\mathbf{r}}
\newcommand{\bbr}{\bar{\mathbf{r}}}

\newcommand{\wL}{\widetilde{L}}
\newcommand{\wsN}{\widetilde{\sN}}
\newcommand{\wrho}{\widetilde{\rho}}

\newcommand{\hy}{\widehat{y}}
\newcommand{\hrho}{\widehat{\rho}}

\newcommand{\Gh}{\widehat{G}}

\newcommand{\CC}{\mathbb{C}}

\newcommand{\NN}{\mathbb{N}}

\newcommand{\QQ}{\mathbb{Q}}

\newcommand{\ZZ}{\mathbb{Z}}

\newcommand{\id}{\textup{id}}
\newcommand{\im}{\textup{im}}

\newcommand{\nin}{\noindent}

\newcommand{\ra}{\rightarrow}
\newcommand{\xra}{\xrightarrow}
\newcommand{\lra}{\longrightarrow}
\newcommand{\co}{\colon\!}

\newcommand{\Gsign}{\textup{G-sign}}
\renewcommand{\mod}{\textup{mod }}
\newcommand{\res}{\textup{res}}
\renewcommand{\max}{\textup{max}}

\newcommand{\RhG}{R_{\widehat G}}
\newcommand{\RhGp}{R_{\widehat G}^+}
\newcommand{\RhGm}{R_{\widehat G}^-}

\newcommand{\pr}{\textup{pr}}

\newcommand{\mnin}{\medskip\noindent}

\title[On Fake lens spaces]{On the classification of fake lens spaces}
\author{Tibor Macko, Christian Wegner}



\address{Mathematisches Institut \\ Universit\"at M\"unster \\
Einsteinstra{\ss}e 62 \\ M\"unster, D-48149 \\ Germany \\ and
Matematick\'y \'Ustav SAV \\ \v Stef\'anikova 49 \\ Bratislava,
SK-81473 \\ Slovakia} \email{macko@uni-muenster.de}
\address{Mathematisches Institut \\ Universit\"at M\"unster \\
Einsteinstra{\ss}e 62 \\ M\"unster, D-48149 \\ Germany}
\email{c.wegner@uni-muenster.de}

\thanks{The authors are supported by SFB 478 Geometrische
Strukturen in der Mathematik, M\"unster.}

\begin{document}

\maketitle

\begin{abstract}
In the first part of the paper we present a classification of fake
lens spaces of dimension $\geq 5$ whose fundamental group is the
cyclic group of order any $N \geq 2$. The classification is stated in
terms of the simple structure sets in the sense of surgery theory.
The results use and extend the results of Wall and others in the
cases $N = 2$ and $N$ odd and the results of the authors of the
present paper in the case $N=2^K$.

In the second part we study the suspension map between the simple
structure sets of lens spaces of different dimensions. As an
application we obtain an inductive geometric description of the
torsion invariants of fake lens spaces.
\end{abstract}

%%%%%%%%%%%%%%%%%%%%%%%%%%%%%%%%%%%%%

\section*{Introduction}

%%%%%%%%%%%%%%%%%%%%%%%%%%%%%%%%%%%%%

A {\it fake lens space} is the orbit space of a free action of a
finite cyclic group $G$ on a sphere $S^{2d-1}$. It is a
generalization of the notion of a {\it lens space} which is the
orbit space of a free action which comes from a unitary
representation. The classification of lens spaces is a classical
topic in algebraic topology and algebraic $K$-theory well explained
for example in \cite{Milnor(1966)}. For the classification of fake
lens spaces in dimension $\geq 5$ methods of surgery theory are
especially suitable. The classification of fake lens spaces with $G$
of order $N = 2$ or $N$ odd was obtained and published in the books
\cite{Wall(1999)}, \cite{LdM(1971)}. In \cite{Macko-Wegner(2008)} we
addressed the problem for $N = 2^K$. In the present paper we treat
the general case $N \geq 2$. The previous classification results are
used as an input.

The general approach to solving the classification problem is the
same in all cases. Firstly, Reidemeister torsion is used to provide
the simple homotopy classification. The homeomorphism classification
within a simple homotopy type can then be formulated in terms of the
{\it simple structure set} $\sS^s(X)$ of a closed $n$-manifold $X$.
An element of $\sS^s(X)$ is represented by a simple homotopy
equivalence $f \co M \ra X$ from a closed $n$-manifold $M$. Two such
$f \co M \ra X$, $f' \co M' \ra X$ are equivalent if there exists a
homeomorphism $h \co M \ra M'$ such that $f' \circ h \simeq f$. The
simple structure set $\sS^s(X)$ is a priori just a pointed set with
the base point $\id \co X \ra X$. However, it can also be endowed
with a preferred structure (in some sense) of an abelian group (see
\cite[chapter 18]{Ranicki(1992)}). Surgery theory is the standard
method for calculating the simple structure set. In the case of fake
lens spaces, in addition to the usual methods of surgery theory, the
so-called reduced $\rho$-invariant and a certain amount of
calculations with this invariant are needed.

We have recalled the above scheme in \cite{Macko-Wegner(2008)} in
detail and we have made the calculations in the case $N = 2^K$. In
the general case we follow the same scheme, but different methods
and calculations are needed in certain stages, in particular those
concerning the reduced $\rho$-invariant. Also the results for $N =
2^K$ and $N$ odd are used as an input. That is the reason why we
present them in a separate paper. Therefore when discussing the
results of the general scheme we will often refer the reader for a
detailed discussion to \cite{Macko-Wegner(2008)} and here we will
concentrate on the calculational issues that are different for
general $N \geq 2$.

The invariants which distinguish the elements of this set and the
homeomorphism classification of polarized fake lens spaces rather
than just elements of the simple structure set are also discussed.
See Introduction and section 1 of \cite{Macko-Wegner(2008)} for a
more detailed general discussion on the relations between these
classification problems.

In the second part of the paper we study the suspension map between
the simple structure sets of lens spaces of different dimensions. As
an application we obtain an inductive geometric description of the
invariants of fake lens spaces on the torsion subgroups of the
simple structure sets calculated in the first part.

\section{Statement of the results}

\begin{defn}
A  {\it fake lens space} $L^{2d-1}(\alpha)$ is a manifold obtained
as the orbit space of a free action $\alpha$ of the group $G =
\ZZ_N$ on $S^{2d-1}$.
\end{defn}

\noindent The fake lens space $L^{2d-1}(\alpha)$ is a
$(2d-1)$-dimensional manifold with $\pi_1 (L^{2d-1}(\alpha)) \cong G
= \ZZ_N$ and universal cover $S^{2d-1}$. The main theorem in this
paper is:

\begin{thm}\label{main-thm}
Let $L^{2d-1}(\alpha)$ be a fake lens space with
$\pi_1(L^{2d-1}(\alpha)) \cong \ZZ_N$ where $N=2^K\cdot M$ with $K
\geq 0$, $M$ odd and $d \geq 3$. Then we have
\[
\sS^s (L^{2d-1}(\alpha)) \cong \bar \Sigma_N (d) \oplus
\bigoplus_{i=1}^{c} \ZZ_{2^{\min\{K,1\}}} \oplus \bigoplus_{i=1}^{c}
\ZZ_{2^{\min\{K,2i\}}}
\]
where $\bar \Sigma_N (d)$ is a free abelian group. If $N$ is odd
then its rank is $(N-1)/2$. If $N$ is even then its rank is $N/2-1$
if $d=2e+1$ and $N/2$ if $d=2e$. In the torsion summand we have $c =
\lfloor (d-1)/2 \rfloor$.
\end{thm}

The isomorphism from the above theorem has an interpretation in
terms of known invariants. The first invariant is the reduced
$\rho$-invariant, denoted $\wrho$, which takes values in $\QQ
\RhG^{(-1)^d}$, the underlying abelian group of the
$(-1)^d$-eigenspace of the rationalized complex representation ring
of $G$ modulo the ideal generated by the regular representation. The
group $\bar \Sigma_N (d)$ is a subgroup of this group. The
invariants $\br_{2i}$ below are closely related to the normal
invariants from surgery theory, see the discussion from
\cite[section 1]{Macko-Wegner(2008)} for more details.

\begin{cor} \label{main-thm-2}
Let $L^{2d-1}(\alpha)$ be a fake lens space with $\pi_1
(L^{2d-1}(\alpha)) \cong \ZZ_N$ where $N = 2^K \cdot M$ with $K \geq
0$, $M$ odd and $d \geq 3$. There exists a collection of invariants
\[
\br_{4i} \co \sS^s (L^{2d-1}(\alpha)) \lra \ZZ_{2^{\min\{K,2i\}}}
\quad \textup{and} \quad \br_{4i-2} \co \sS^s (L^{2d-1}(\alpha))
\lra \ZZ_{2^{\min\{K,1\}}}
\]
where $1 \leq i \leq c = \lfloor (d-1)/2 \rfloor$ which together
with the $\wrho$-invariant induce a one-to-one correspondence
between elements $a \in \sS^s(L^{2d-1}(\alpha))$ and
\begin{enumerate}
\item $\wrho (a) \in \bar \Sigma \subset \QQ \RhG^{(-1)^d}$
\item $\br_{2i} (a) \in \ZZ_{2^{\min\{K,1\}}}, \ZZ_{2^{\min\{K,2i\}}}$.
\end{enumerate}
\end{cor}

To obtain classification of fake lens spaces rather than
classification of elements of the simple structure set we follow
Wall \cite[chapter 14E]{Wall(1999)} and work with polarized fake
lens spaces. These are fake lens spaces equipped with a choice of
orientation and a choice of a generator of the fundamental group,
see \cite[Definition 2.2]{Macko-Wegner(2008)}. The simple homotopy
type of a polarized fake lens space $L^{2d-1} (\alpha)$ is given uniquely by its Reidemeister torsion $\Delta (L^{2d-1} (\alpha))$, which is a unit in $\QQ R_G$, the rational group ring of $G$ modulo the ideal generated by the norm element, as described in \cite[Proposition 2.3]{Macko-Wegner(2008)}. That means that for two polarized fake lens spaces $L^{2d-1} (\alpha)$ and $L^{2d-1} (\beta)$ with the fundamental group $G$ there is a simple homotopy equivalence
$f_{\alpha,\beta} \co L^{2d-1}(\alpha) \ra L^{2d-1} (\beta)$ of
polarized fake lens spaces unique up to homotopy if and only if the
Reidemeister torsions of $L^{2d-1}(\alpha)$ and $L^{2d-1}(\beta)$
coincide. The simple homotopy equivalence $f_{\alpha,\beta}$ gives
us an element of the simple structure set $\sS^s (L^{2d-1}(\beta))$.
We can formulate the classification as follows:

\begin{cor} \label{main-thm-3}
Let $L^{2d-1} (\alpha)$ and $L^{2d-1} (\beta)$ be polarized lens
spaces with the fundamental group $G = \ZZ_N$, where $N = 2^K \cdot
M$ with $K \geq 0$, $M$ odd and $d \geq 3$. There exists a polarized
homeomorphism between $L^{2d-1} (\alpha)$ and $L^{2d-1} (\beta)$ if
and only if
\begin{enumerate}
 \item $\Delta (L^{2d-1}(\alpha)) = \Delta (L^{2d-1} (\beta))$,
 \item $\rho (L^{2d-1}(\alpha)) = \rho (L^{2d-1} (\beta))$,
 \item $\br_{2i} (f_{\alpha,\beta}) = 0$ for all $i$.
\end{enumerate}
\end{cor}

In the second part of the paper we study the so-called suspension
homomorphisms $\Sigma \co \sS^s (L^{2d-1}(\alpha_k)) \ra
\sS^s(L^{2d+1}(\alpha_k))$, see \cite[section
3.3]{Macko-Wegner(2008)}. Here $\alpha_k$ denotes certain actions of
$\ZZ_N$ on $S^{2d-1}$ and on $S^{2d+1}$which yield standard lens spaces, see definition in \cite[section 2]{Macko-Wegner(2008)}. The reason for restricting to the case $\alpha = \alpha_k$ is explained in Remark \ref{restriction-to-lens-spaces}. The main results are Theorem
\ref{thm1} and \ref{thm2}. Their statements require certain notation
and are therefore relegated to section \ref{sec:suspension}.
Nevertheless, we can now state a corollary of these two theorems.

\begin{cor}
Let $N=2^K\cdot M$ with $K \geq 1$, $M$ odd. If $e \geq 1$ then
there is an exact sequence:
\[
0 \ra \sS^s (L^{4e+1} (\alpha_k)) \xra{\Sigma} \sS^s (L^{4e+3}
(\alpha_k)) \ra \ZZ \ra 0.
\]
If $e \geq 2$ then there is an exact sequence
\[
0 \ra \ZZ \ra \sS^s (L^{4e-1} (\alpha_k)) \xra{\Sigma} \sS^s
(L^{4e+1} (\alpha_k)) \ra \ZZ_2 \ra 0.
\]
\end{cor}

As an application we obtain an improvement on the invariants from
Corollary \ref{main-thm-3}. Although the invariants $\br_{4i-2}$
have certain geometric description, which we recall in Remark
\ref{rem:splitting-invariants-4i-2}, we are not able to give such a
description for the invariants $\br_{4i}$. However, using our
results about the suspension map we obtain new invariants
$\bbr_{4i}$. Although, as explained in section \ref{sec:invariants},
there is still a certain choice involved, the invariants $\bbr_{4i}$
are an improvement, since after that choice is made, they have an
inductive geometric interpretation. A more detailed statement again
needs some notation and is therefore given in section
\ref{sec:invariants} as Corollary \ref{cor:splitting-invariants},
and Remark \ref{rem:splitting-invariants-4i}.

The paper is organized as follows. In section \ref{sec:ses} we
recall the information about the simple structure set of fake lens
spaces which we obtain from surgery theory. In section
\ref{sec:rho-invariant} we recall the definition of the
$\wrho$-invariant and we prove some of its properties, in particular
we show that it gives a group homomorphism of the simple structure
set. In section \ref{sec:revision} we recall the results in the
cases $N = 2^K$ and $N = M$ odd which are used in section
\ref{sec:calculations}, which contains the calculations. In section
\ref{sec:suspension} the suspension map is studied and in the final
section \ref{sec:invariants} we provide the geometric description of
the torsion invariants.

\section{The surgery exact sequence} \label{sec:ses}

%%%%%%%%%%%%%%%%%%%%%%%%%%%%%%%%%%%%%%%%%

For the simple homotopy classification of fake lens spaces we refer
the reader to \cite[section 2]{Macko-Wegner(2008)} where we review
the results of Wall. Here we concentrate on the homeomorphism
classification within a simple homotopy type. The main tool is the
surgery exact sequence computing the structure set $\sS^s(X)$ for a
given $n$-manifold $X$ with $n \geq 5$:
\begin{equation} \label{ses}
\cdots \ra \sN_\partial (X \times I) \xra{\theta} L^s_{n+1} (G)
\xra{\partial} \sS^s (X) \xra{\eta} \sN(X) \xra{\theta} L^s_n (G),
\end{equation}
where $G = \pi_1 (X)$. The terms in the sequence are reviewed in the
detail in \cite[section 3]{Macko-Wegner(2008)}. Now we analyze them
for $X = L^{2d-1}_N(\alpha)$. The following proposition from
\cite{Hambleton-Taylor(2000)} describes the $L$-theory needed. Here
$R_{\CC} (G)$ denotes the complex representation ring of a group $G$
and the superscripts $\pm$ denote the $\pm$-eigenspaces with respect
to the involution given by complex conjugation.

\begin{thm} \label{L(G)}
For $G = \ZZ_N$ we have that
\begin{align*}
L^s_n (G) & \cong
\begin{cases}
4 \cdot R_{\CC}^+ (G) & n \equiv 0 \; (\mod 4) \; (\Gsign, \;
\mathrm{purely} \; \mathrm{real}) \\
0 & n \equiv 1 \; (\mod 4) \\
4 \cdot R_{\CC}^- (G) \oplus \ZZ_2 & n \equiv 2 \; (\mod 4) \;
(\Gsign, \; \mathrm{purely}
\; \mathrm{imaginary}, \mathrm{Arf}) \\
\ZZ_2 & n \equiv 3 \; (\mod 4) \; (\mathrm{codimension} \; 1 \;
\mathrm{Arf})
\end{cases} \\
\widetilde L^s_{2k} (G) & \cong 4 \cdot \RhG^{(-1)^k} \;
\textit{where} \; \RhG^{(-1)^k} \; \textit{is} \; R_{\CC}^{(-1)^k}
(G) \; \textit{modulo the regular representation.}
\end{align*}
\end{thm}

For the normal invariants we have in general

\begin{cor} \label{ni}
If $X$ is rationally trivial we have an isomorphism of abelian
groups
\begin{align*}
\sN(X) & \cong \sN(X)_{(2)} \oplus \sN(X)_{(odd)} \\
& \cong \Big( \bigoplus_{i \geq 1} \big( H^{4i} (X ; \ZZ_{(2)})
\oplus H^{4i-2} (X;\ZZ_2) \big) \Big) \oplus KO (X) \otimes
\ZZ[\frac{1}{2}]
\end{align*}
\end{cor}

\nin When $X$ is a fake lens space $L^{2d-1}(\alpha)$ with $\pi_1
(L^{2d-1}(\alpha)) \cong G = \ZZ_N$ we obtain
\begin{align*} \label{ni-lens-spaces}
\sN (L^{2d-1}(\alpha)) \cong & \bigoplus_{i=1}^{\lfloor (d-1)/2
\rfloor } H^{4i} (L^{2d-1}(\alpha);\ZZ_{(2)}) \oplus
\bigoplus_{i=1}^{\lfloor d/2 \rfloor} H^{4i-2}
(L^{2d-1}(\alpha);\ZZ_2) \\
& \oplus KO (L^{2d-1}(\alpha)) \otimes \ZZ\Big[\frac{1}{2}\Big]
\end{align*}
The first two summands can be easily calculated explicitly and we
denote the factors
\begin{align}
\bt_{4i} & \co \sN(L^{2d-1}(\alpha)) \ra H^{4i} (L^{2d-1}(\alpha);\ZZ_{(2)}) \cong \ZZ_{2^K} \\
\bt_{4i-2} & \co \sN(L^{2d-1}(\alpha)) \ra H^{4i-2}
(L^{2d-1}(\alpha);\ZZ_2) \cong \ZZ_2.
\end{align}
The last summand is more difficult to calculate explicitly, but we
will not need the exact calculation. Note that when Wall analyzed
the case $N = M$ he only needed that the order of the group $KO
(L^{2d-1}(\alpha)) \otimes \ZZ[\frac{1}{2}]$ is $M^c$ with $c =
\lfloor (d-1)/2 \rfloor$, which is an easy Atiyah-Hirzebruch
spectral sequence argument. This will also be sufficient for us. We
will still need some more notation, so the projection onto this last
summand will be denoted
\begin{equation}
\bt_{(odd)} \co \sN(L^{2d-1}(\alpha)) \ra KO (L^{2d-1}(\alpha))
\otimes \ZZ\Big[\frac{1}{2}\Big].
\end{equation}
We will also sometimes put together the $2$-local invariants and
denote $\bt_{(2)} = (\bt_{2i})_i$ and finally $\bt =
(\bt_{(2)},\bt_{(odd)})$. These projections will also be sometimes
used to identify the elements of $\sN(L^{2d-1}(\alpha))$ as $t =
((t_{2i})_i,t_{(odd)})$. Even more information is obtained from the
following
\begin{thm}[\cite{Wall(1999)}]
\
\begin{enumerate}
\item If $d=2e$ then the map
\[
\theta \co \sN(L^{2d-1}(\alpha)) \ra L^s_{2d-1}(G) = L^s_{4e-1}(G) =
\ZZ_2
\]
is given by $\theta (x) = \bt_{4e-2} (x) \in  \ZZ_2$.
\item
The map
\[
\theta \co \sN_\partial(L^{2d-1}(\alpha) \times I) \ra L^s_{2d}(G)
\]
maps onto the summand $L^s_{2d}(1)$.
\end{enumerate}
\end{thm}
\nin  Hence we obtain the short exact sequence
\begin{equation} \label{ses-lens-2d-1}
0 \ra \wL^s_{2d} (G) \xra{\partial} \sS^s   (L^{2d-1}(\alpha))
\xra{\eta} \widetilde{\sN}(L^{2d-1}(\alpha)) \ra 0
\end{equation}
where
\begin{align*}
\widetilde{\sN}(L^{4e-1}(\alpha)) & = \mathrm{ker} \; \big (
\bt_{4e-2} \co
{\sN}(L^{4e-1}(\alpha)) \ra H^{4e-2} (L^{4e-1}(\alpha);\ZZ_2) \cong \ZZ_2 \big ), \\
\widetilde{\sN}(L^{4e+1}(\alpha)) & = \sN(L^{4e+1}(\alpha)).
\end{align*}
in other words
\begin{equation} \label{red-ni-lens-spaces}
\widetilde \sN (L^{2d-1}(\alpha)) \cong \bigoplus_{i=1}^{c} \ZZ_{2^K}
\oplus \bigoplus_{i=1}^{c} \ZZ_2 \oplus KO (L^{2d-1}(\alpha))
\otimes \ZZ\Big[\frac{1}{2}\Big]
\end{equation}
where $c = \lfloor (d-1)/2 \rfloor$ and where the order of the last
summand is $M^c$. The first term in the sequence
(\ref{ses-lens-2d-1}) is understood by Theorem \ref{L(G)}, the third
term is understood by (\ref{red-ni-lens-spaces}). Hence we are left
with an extension problem.

\begin{rem}
We will also need to work with fake complex projective spaces.
Wall's calculation of $\sS (\CC P^{d-1})$ and $\sN (\CC P^{d-1})$ is
reviewed in \cite[subsection 3.1]{Macko-Wegner(2008)}. Recall that
$\sS (\CC P^{d-1})$ is a subgroup of $\sN (\CC P^{d-1})$. Both
groups are calculated in terms of invariants $\bs_{4i} \in \ZZ$ and
$\bs_{4i-2} \in \ZZ_2$.
\end{rem}

In this paper we will also make systematic use of the functoriality
of the normal invariants and of the structure set. Let $H < G \leq
S^1$ be an inclusion of subgroups and let $\alpha$ be a free action
of $G$ on $S^{2d-1}$. The inclusion $H < G$ induces a free action of
$H$ which we also denote $\alpha$. Let now $L_G(\alpha)$ denote the
resulting $(2d-1)$-dimensional fake lens spaces in case $G < S^1$ or
the resulting $(2d-2)$-dimensional fake complex projective space in
case $G = S^1$. We have fiber bundles
\[
 p_H^G \co L_H(\alpha) \lra L_G(\alpha).
\]
The maps $p_H^G$ induce the vertical maps in the following diagram
\[
\xymatrix{
 \sS^s (L_G(\alpha)) \ar[r]^{\eta} \ar[d]_{(p_H^G)^{!}} & \sN (L_G(\alpha)) \ar[d]^{(p_H^G)^{!}} \\
 \sS^s (L_H(\alpha)) \ar[r]^{\eta} & \sN (L_H(\alpha))
}
\]
The right hand map preserves the localization at $2$ and away from
$2$ and we have that $(p_H^G)^{!} \co \sN (L_G(\alpha))_{(2)} \lra
\sN (L_H(\alpha))_{(2)}$ is given by reduction modulo $|H|$. Suppose
now that $G = \ZZ_N$, $N = 2^K \cdot M$, $M$ odd, and $H =
\ZZ_{N'}$, $N' = 2^{K'} \cdot M'$, $M'$ odd. If $M = M'$ then
$(p_H^G)^{!} \co \sN (L_G(\alpha))_{(odd)} \lra \sN
(L_H(\alpha))_{(odd)}$ is an isomorphism. If $K = K'$ then
$(p_H^G)^{!} \co \sN (L_G(\alpha))_{(2)} \lra \sN
(L_H(\alpha))_{(2)}$ is an isomorphism. In fact we have for $N = 2^K
\cdot M$ with $M$ odd that
\begin{equation} \label{ni-N-vs-2K-M}
 p_{2^K}^N \oplus p_{M}^N \co \sN (L_N) \cong \sN (L_N)_{(2)} \oplus \sN (L_N)_{(odd)} \lra \sN (L_{2^K}) \oplus \sN (L_M)
\end{equation}
is an isomorphism. (We have left out $(\alpha)$ everywhere due to
the lack of space.) In the sequel we will sometimes use the notation $L_G (\alpha)$ when we feel the need for specifying $G$. In the other cases the group involved should be clear.

Recall from \cite[section 2]{Macko-Wegner(2008)} the definition of
the actions $\alpha_k$ of $\ZZ_N$ on $S^{2d-1}$ for $k \in \NN$,
$(k,N) =1$, which yield standard lens spaces $L^{2d-1}(\alpha_k)$.
Further recall from \cite[subsection 3.3]{Macko-Wegner(2008)} the
join of fake lens spaces. The join with $L^1(\alpha_k)$ defines a
homomorphism $\Sigma_k \co \sS^s (L^{2d-1}(\alpha_1)) \lra
\sS^s(L^{2d+1}(\alpha_k))$. The inclusion $L^{2d-1}(\alpha_1)
\subset L^{2d+1}(\alpha_k)$ induces a restriction map on the groups
of normal invariants denoted by $\res \co \sN (L^{2d+1}(\alpha_k))
\lra \sN (L^{2d-1}(\alpha_1))$ and we have a commutative diagram
\cite[Lemma 14A.3]{Wall(1999)}:
\begin{equation} \label{susp-diagram}
\begin{split}
\xymatrix{
\sS^s (L^{2d-1}(\alpha_1)) \ar[r]^{\eta} \ar[d]_{\Sigma_k} & \sN (L^{2d-1}(\alpha_1)) \\
\sS^s (L^{2d+1}(\alpha_k)) \ar[r]^{\eta} & \sN (L^{2d+1}(\alpha_k))
\ar[u]_{\res} }
\end{split}
\end{equation}
Note that we have $t_{2i} = \res (t_{2i})$. Moreover, the map
\begin{equation}
\res \co \widetilde \sN(L^{2d+1}(\alpha_1)) \lra \widetilde
\sN(L^{2d-1} (\alpha_1))
\end{equation}
is an isomorphism when $d = 2e+1$ and it is onto when $d = 2e$. A
similar diagram exists for the situation $\CC P^d = \CC P^{d-1} \ast
\mathrm{pt}$.

%%%%%%%%%%%%%%%%%%%%%%%%%%%%%%%%%%%%%%%%%

\section{The $\rho$-invariant} \label{sec:rho-invariant}

%%%%%%%%%%%%%%%%%%%%%%%%%%%%%%%%%%%%%%%%%

Similarly as in \cite{Macko-Wegner(2008)} we will use the
$\rho$-invariant of odd-dimensional manifolds to solve our extension
problem. The definition we use is the same as in
\cite{Macko-Wegner(2008)}, but the formulation of the properties and
their proofs have to be adjusted to the more general case. This is
the content of the present section.

%%%%%%%%%%%%%%%%%%%%%%%%%%%%%%%%%%%%%%%%%

\subsection{Definitions}

%%%%%%%%%%%%%%%%%%%%%%%%%%%%%%%%%%%%%%%%%

\

\  \noindent Recall the definition of the $\rho$-invariant (see also
\cite[subsection 4.1]{Macko-Wegner(2008)}).

\begin{defn}{\cite[Remark after Corollary 7.5]{Atiyah-Singer-III(1968)}} \label{defn-rho-1}
Let $X^{2d-1}$ be a closed manifold with $\pi_1 (X) \cong G$ a
finite group. Define
\begin{equation}
\rho (X) = \frac{1}{k} \cdot \Gsign(\widetilde Y) \in \QQ R^{(-1)^d}
(G)/ \langle \textup{reg} \rangle
\end{equation}
for some $k \in \NN$ and $(Y,\partial Y)$ such that $\pi_1 (Y) \cong
\pi_1 (X)$ and $\partial Y = k \cdot X$. The symbol $\langle
\textup{reg} \rangle$ denotes the ideal generated by the regular
representation.
\end{defn}
We remind the reader that there is also another definition which
works for actions of compact Lie groups, in particular for
$S^1$-actions, on certain odd-dimensional manifolds. Whenever the
two definitions apply, they coincide. For $G < S^1$ we identify
$R(G)$ with $\ZZ \Gh$ and we adopt the notation $\RhG :=  R(G) /
\langle \textup{reg} \rangle$ and $\RhG = \ZZ [\chi] / \langle 1 +
\chi + \cdots + \chi^{N-1} \rangle$ as explained in \cite[section
4.1]{Macko-Wegner(2008)}.

Hence we have the $\rho$-invariant defined for fake lens spaces and
for fake complex projective spaces. We continue with a list of some
basic properties. For the join $L \ast L'$ of fake lens spaces $L$
and $L'$ we have \cite[chapter 14A]{Wall(1999)}
\begin{equation} \label{rho-join}
\rho (L \ast L') = \rho (L) \cdot \rho (L').
\end{equation}
For $L^1(\alpha_k)$ we have \cite[Proof of Theorem
14C.4]{Wall(1999)}
\begin{equation} \label{rho-alpha-k}
\rho (L^1(\alpha_k)) = f_k \in \QQ \RhGm
\end{equation}
where $f_k$ is defined as follows.
\begin{defn}
For $k \in \NN$ with $(N,k)=1$ we set
\[
f_k := \frac{1+\chi^k}{1-\chi^k}
\]
and
\begin{equation*}
\everymath{\displaystyle} f'_k := \left\{
\begin{array}{ll}
 \frac{1-\chi+\chi^2-\cdots-\chi^{k-2}+\chi^{k-1}}{1+\chi+\chi^2+\cdots+\chi^{k-2}+\chi^{k-1}} & \mathrm{\; for \;} k \mathrm{\; odd}, \\
 \rule{0in}{5ex}
 \frac{\chi^k-\chi^{k+1}+\chi^{k+2}-\cdots-\chi^{N-2}+\chi^{N-1}}{1+\chi+\chi^2+\cdots+\chi^{k-2}+\chi^{k-1}} & \mathrm{\; for \;} k \mathrm{\; even}.
\end{array}
\right.
\end{equation*}
We abbreviate $f := f_1$. Note that if $k$ is even, we necessarily
have $N = M$ odd.
\end{defn}

\begin{lem}\label{f_k-lem}
Let $G = \ZZ_N$ with $N=2^K \cdot M$, $M$ odd. For $k \in \NN$ with
$(N,k)=1$ we have
\[
f_k \in \QQ \RhGm, \qquad f_k = f \cdot f'_k, \qquad f'_k \in \RhG.
\]
\end{lem}
\begin{proof}
Notice that $1-\chi^k$ is invertible in $\QQ \RhG$ because
\[
(1-\chi^k)^{-1} = -\frac{1}{N}(1 + 2 \cdot \chi^k + 3 \cdot
\chi^{2k} + \cdots + N \cdot \chi^{(N-1)k}) \in \QQ \RhG.
\]
Therefore, $f_k \in \QQ \RhG$ and the identity
\[
\frac{1+\chi^{-k}}{1-\chi^{-k}} = - \frac{1+\chi^k}{1-\chi^k} = -
f_k
\]
implies $f_k \in \QQ \RhGm$. An easy calculation shows $f_k = f
\cdot f'_k$. That $f'_k \in \RhG$ follows from the fact that
$1+\chi+\chi^2+\cdots+\chi^{k-1}$ is invertible in $\RhG$. The
inverse is given by $1+\chi^k+\chi^{2k}+\cdots+\chi^{(r-1)k}$ where
$r$ denotes a natural number such that $r \cdot k - 1$ is a multiple
of $N$.
\end{proof}

\begin{rem} \label{chi-rem}
Recall that for $G < S^1$ we have a canonical isomorphism $R(G) =
\ZZ \widehat G$. Suppose that we also have a subgroup $H$ of $G$ and
denote the inclusion $i \co H \hookrightarrow G$. Then the ring
homomorphism $R(G) \ra R(H)$ induced by the restriction is
identified with the ring homomorphism $\ZZ \widehat G \ra \ZZ
\widehat H$ induced by the group homomorphism $\widehat i \co
\widehat G \ra \widehat H$. This homomorphism sends a generator to a
generator. The choice of a generator $\chi$ of $\widehat G$ also
gives us an identification $R(G) = \ZZ \widehat G =  \ZZ[\chi] /
\langle \chi^{|G|} - 1\rangle$. Hence we can think of the induced
homomorphism $R(G) \lra R(H)$ as of the obvious quotient map
$\ZZ[\chi] / \langle \chi^{|G|} - 1\rangle \lra \ZZ[\chi] / \langle
\chi^{|H|} - 1\rangle$.
\end{rem}

%%%%%%%%%%%%%%%%%%%%%%%%%%%%%%%%%%%%%%%%%

\subsection{Homomorphism}

\

%%%%%%%%%%%%%%%%%%%%%%%%%%%%%%%%%%%%%%%%%

\

The $\rho$-invariant defines a function of $\sS^s (X)$ by sending $a
= [h \co M \lra X]$ to $\wrho (a) = \rho(M) - \rho (X)$. If we put
on $\sS^s (X)$ the abelian group structure from \cite[chapter
18]{Ranicki(1992)} it is not clear whether $\wrho$ is a homomorphism
in general. It is the aim of the present subsection to show this
claim for $X = L^{2d-1}(\alpha)$. More precisely we prove the
following

\begin{prop}\label{ses-vs-rep-thy}
There is the following commutative diagram of abelian groups and
homomorphisms with exact rows
\[
\xymatrix{ 0 \ar[r] & \wL^s_{2d} (G) \ar[r]^(0.4){\partial}
\ar[d]_{\cong}^{\Gsign} & \sS^s (L^{2d-1}(\alpha)) \ar[r]^{\eta}
\ar[d]^{\widetilde \rho}&
\widetilde \sN (L^{2d-1}(\alpha)) \ar[r] \ar[d]^{[\widetilde \rho]}& 0 \\
0 \ar[r] & 4 \cdot R^{(-1)^d}_{\widehat G} \ar[r] & \QQ
R^{(-1)^d}_{\widehat G} \ar[r] & \QQ R^{(-1)^d}_{\widehat G}/ 4
\cdot R^{(-1)^d}_{\widehat G} \ar[r] & 0 }
\]
where $[\wrho]$ is the homomorphism induced by $\wrho$.
\end{prop}

Similarly as in \cite{Macko-Wegner(2008)} the commutativity of the
left hand square follows from \cite[Theorem 2.3]{Petrie(1970)}. We
need to show that $\wrho$ and  $[\wrho]$ are homomorphisms. This
will be proved first for $\alpha_1$, then for $\alpha_k$, and
finally for general $\alpha$.

We start with a useful lemma.

\begin{lem} \label{f-inverse}
Let $N = 2^K \cdot M$, $M$ odd, $G = \ZZ_N$. Then there exists $g
\in \QQ \RhG^-$ such that for all $x \in \QQ \RhG^-$ we have $x = g
\cdot f \cdot x$.
\end{lem}

\begin{proof}
We use the Chinese remainder theorem which tell us that we have an
isomorphism of rings
\begin{multline}
\QQ [\chi] / \langle 1 + \chi + \cdots + \chi^{N-1} \rangle
\xra{\cong} \\ \xra{\cong} \bigoplus_{l = 0}^{K-1} \QQ [\chi] /
\langle 1 + \chi^{2^l} \rangle \oplus \QQ [\chi] / \langle 1 +
(\chi^{2^K}) + \cdots + (\chi^{2^K})^{M-1} \rangle
\end{multline}
Denote the projections $\pr_l \co \QQ [\chi] / \langle 1 + \chi +
\cdots + \chi^{N-1} \rangle \lra \QQ [\chi] / \langle 1 + \chi^{2^l}
\rangle$ and $\pr \co \QQ [\chi] / \langle 1 + \chi + \cdots +
\chi^{N-1} \rangle \lra \QQ [\chi] / \langle 1 + (\chi^{2^K}) +
\cdots + (\chi^{2^K})^{M-1} \rangle$.

\medskip\nin If $K=0$ then $f$ is invertible since
\[
(1+\chi)^{-1} = (1 + \chi^2 + \cdots + \chi^{M-2}) \in \ZZ[\chi] /
\langle 1 + \chi + \cdots \chi^{N-1} \rangle.
\]

\medskip\nin If $K > 0$ and $l > 0$ then notice that $(1 + \chi)$ is invertible in $\QQ
[\chi] / \langle 1 + \chi^{2^l} \rangle$:
\[
h_l : = (1+\chi)^{-1} = \frac{1}{2} \cdot A_l \textrm{\; with \;}
A_l := 1-\chi+\chi^2-\chi^3+\cdots -\chi^{2^l-1}.
\]
It is also invertible in $\QQ [\chi] / \langle 1 + (\chi^{2^K}) +
\cdots + (\chi^{2^K})^{M-1} \rangle$:
\[
h : = (1+\chi)^{-1} = -\frac{A_K}{M} \cdot \big(1 - 2 \cdot
(\chi^{2^K}) + 3 \cdot (\chi^{2^K})^2 - \cdots + M \cdot
(\chi^{2^K})^{M-1} \big).
\]
Let $g \in \QQ \RhG^-$ be such that $\pr_l (g) = h_l$ for $l > 0$,
$\pr (g) = h$ and $\pr_0 (g) = 0$. Further recall that we can write
$x \in \QQ \RhG^-$ as
\[
x = \sum_{r=1}^{N/2-1} a_r \cdot (\chi^r - \chi^{N-r}).
\]
with $a_r \in \QQ$. Since $\chi^r - \chi^{N-r}$ is a multiple of $1
+ \chi$, we conclude $\pr_0 (x) = 0$. We obtain $\pr_l (gfx) = \pr_l
(x)$ for $l >0$, $\pr(gfx)= \pr(x)$ and $\pr_0 (gfx) = 0 =
\pr_0(x)$. Hence $gfx = x$.
\end{proof}

\begin{prop} \label{ni-wrho-homomorphism}
The function $[\wrho] \co \sN (L_N^{2d-1} (\alpha_1)) \lra
\QQ\RhG^{(-1)^d}/4 \cdot \RhG^{(-1)^d}$ is a homomorphism.
\end{prop}

\begin{proof}
\mnin Consider the following commutative diagram: %\ar@/_4pc/ [dd]
\begin{equation} \label{res-diag}
\begin{split}
\xymatrix{
0 \ar[r] & \sS (\CC P^{d-1} ) \ar[d]^{p^{!}} \ar[r]^{\eta} & \sN (\CC P^{d-1}) \ar[d]^{p^{!}} \ar[r] & L_{2(d-1)} (1) \\
\wL^s_{2d} (G) \ar[r] \ar[d]^{\Gsign} & \sS^s (L^{2d-1}(\alpha_1))
\ar[r]^{\eta}
\ar[d]^{\wrho} & \widetilde \sN (L^{2d-1}(\alpha_1)) \ar[r] \ar[d]^{[\wrho]} & 0  \\
4 \cdot \RhG^{(-1)^d} \ar[r] & \QQ \RhG^{(-1)^d} \ar[r] & \QQ
\RhG^{(-1)^d}/4 \cdot \RhG^{(-1)^d} }
\end{split}
\end{equation}
where the maps $p^{!}$, $\eta$ are homomorphisms. The composition
$\wrho \circ p^{!}$ is a homomorphism by \cite[Theorem
14C.4]{Wall(1999)}. Hence also the composition $[\wrho] \circ \eta
\circ p^{!}$ is a homomorphism.

If $d=2e$ then $\eta \circ p^{!}$ is surjective which implies that
$[\wrho]$ is a homomorphism. If $d = 2e+1$ then $\eta \circ p^{!}$
is not surjective. But we have the following commutative diagram:
\begin{equation}
\begin{split}
\xymatrix{
\sS^s (L^{4e+1}(\alpha_1)) \ar[r]^{\eta} \ar[d]_{\Sigma_1} & \widetilde{\sN} (L^{4e+1}(\alpha_1)) \\
\sS^s (L^{4e+3}(\alpha_1)) \ar[r]^{\eta} & \widetilde{\sN}
(L^{4e+3}(\alpha_1)) \ar[u]_{\res} }
\end{split}
\end{equation}
where $\res$ is bijective and we have $\wrho (\Sigma_1 (x)) = f
\cdot \wrho (x)$. Now we can use Lemma \ref{f-inverse} and
calculate:
\begin{align*}
[\wrho] (\eta(x)+\eta(y)) & = [\wrho] (\eta(x+y)) = [\wrho (x+y)] = [g \cdot f \cdot \wrho (x+y)] \\
 &= [g \cdot \wrho (\Sigma_1(x+y))] = [g \cdot \wrho (\Sigma_1 (x) + \Sigma_1 (y))] \\
 &=  [g] \cdot [\wrho] (\eta \Sigma_1 (x) + \eta \Sigma_1 (y)) = [g] \cdot [\wrho] (\eta \Sigma_1 (x)) + [\wrho] (\eta \Sigma_1 (y)) \\
 & = [g] \cdot ([\wrho (\Sigma_1 (x))] + [\wrho (\Sigma_1 (y))]) = [g] \cdot [f] \cdot ([\wrho (x)] + [\wrho (y)]) \\
 & = [\wrho] (\eta (x)) + [\wrho] (\eta (y)).
\end{align*}
This finishes the proof.
\end{proof}

\begin{cor} \label{rho-homomorphism}
The function $\wrho \co \sS^s(L^{2d-1}(\alpha_1)) \lra \QQ
\RhG^{(-1)^d}$ is a homomorphism.
\end{cor}

\begin{proof}
It is enough to show that for every $y$, $y' \in \widetilde \sN
(L^{2d-1}(\alpha_1))$ there exist elements (not necessarily unique)
$a$, $a'$ in $\sS^s (L^{2d-1}(\alpha_1))$ such that $\eta(a) = y$,
$\eta(a') = y'$ and $\wrho (a+a') = \wrho (a) + \wrho (a')$. If this
holds then for any $x,x' \in \sS^s (L^{2d-1}(\alpha_1))$ choose $a$
and $a'$ as above corresponding to the classes $\eta(x)$, $\eta(x')
\in \sN (L^{2d-1} (\alpha_1))$. Then $x = a +\partial(b)$ and $x' =
a' +
\partial(b')$ for some $b$, $b' \in \wL^s_{2d} (G)$ and
\begin{align*}
\wrho (x+x') = & \wrho (a + \partial (b) + a' + \partial (b')) =
\wrho (a + a') +
\wrho (\partial (b) + \partial (b')) \\
= &  \wrho (a) + \wrho(a') + \wrho (\partial (b)) + \wrho(\partial
(b')) = \wrho (x) + \wrho(x').
\end{align*}

When $d=2e$ we can associate to a given $y \in \sN
(L^{2d-1}(\alpha_1))$ an $a \in \sS^s (L^{2d-1}(\alpha_1))$ coming
from the $\sS (\CC P^{d-1})$, i.e. $a = p^{!} (b)$ where $b \in \sS
(\CC P^{d-1})$ such that $p^{!} (\eta (b)) = y$. When we have $y$,
$y'\in \sN (L^{2d-1}(\alpha_1))$ then $\wrho (a+a') = \wrho (p^{!}
(b) + p^{!} (b')) = \wrho (p^{!} (b + b')) = \res (\wrho_{S^1} (b +
b')) = \res (\wrho_{S^1} (b) + \wrho_{S^1} (b')) = \res (\wrho_{S^1}
(b)) + \res(\wrho_{S^1} (b')) = \wrho (a) + \wrho (a')$. Here $\res$
denotes the map on the representation rings induced by the inclusion
$G < S^1$.

When $d = 2e+1$ it follows from Proposition
\ref{ni-wrho-homomorphism} that for $a, a' \in \sS^s
(L^{2d-1}(\alpha_1))$ we have $\wrho (a+a') = \wrho (a) + \wrho (a')
+ z$ for some $z \in 4 \cdot \RhG^-$. Our task is to show $z=0$ for
any choice of $a$, $a'$. We use the fact that $\Sigma_1$ is a
homomorphism and that we have already proved the claim for $d =
2e+2$. That implies $\wrho (\Sigma_1 (a+a')) = \wrho (\Sigma_1 (a) +
\Sigma_1 (a')) = \wrho (\Sigma_1 (a)) + \wrho (\Sigma_1 (a')) = f
\cdot \wrho (a) + f \cdot \wrho(a')$. On the other hand $\wrho
(\Sigma_1 (a+a')) = f \cdot \wrho (a+a') = f \cdot \wrho (a) + f
\cdot \wrho(a') + f \cdot z$. Hence $f \cdot z = 0$ and therefore $z
= g \cdot f \cdot z = 0$ by Lemma \ref{f-inverse}.
\end{proof}

Now we proceed to the case of $\alpha_k$ where $k \in \NN$ is such
that $(N,k)=1$.

\begin{prop}
The function $[\wrho] \co \sN (L^{2d-1}(\alpha_k)) \lra \QQ
\RhG^{(-1)^d}/4 \cdot \RhG^{(-1)^d}$ is a homomorphism.
\end{prop}

\begin{proof}
We will use the result for $\alpha_1$ and the homeomorphisms
\[
L^{2d+1}(\alpha_k) \cong L^{2d-1}(\alpha_1) \ast L^1(\alpha_k) \quad
\mathrm{and} \quad L^{2d+1}(\alpha_k) \cong L^{2d-1}(\alpha_k) \ast
L^1(\alpha_1)
\]
and formulas (\ref{rho-join}) and (\ref{rho-alpha-k}). For $d=2e$
recall the diagram
\begin{equation} \label{join-1-k}
\begin{split}
\xymatrix{ \QQ \RhG^- \ar[d]^{\cdot f_k} & \sS^s (L^{4e-3}
(\alpha_1)) \ar[l]_(0.6){\wrho}
\ar[d]^{\Sigma_k} \ar@{-{>>}}[r]^{\eta} & \widetilde \sN (L^{4e-3} (\alpha_1)) \\
\QQ \RhG^+ & \sS^s (L^{4e-1} (\alpha_k)) \ar[l]_(0.6){\wrho}
\ar@{-{>>}}[r]^{\eta} & \widetilde \sN (L^{4e-1} (\alpha_k))
\ar[u]_{\res}^{\cong} }
\end{split}
\end{equation}
It follows that the composition $[\wrho] \circ (\res)^{-1} \circ
\eta$ is a homomorphism. Since the composition $(\res)^{-1} \circ
\eta$ is surjective, we obtain that the function $[\wrho] \co \sN
(L^{4e-1} (\alpha_k)) \ra \QQ \RhG^{(-1)^d}/4 \cdot \RhG^{(-1)^d}$
is a homomorphism. For $d=2e+1$ recall the diagram
\begin{equation} \label{join-k-k}
\begin{split}
\xymatrix{ \QQ \RhG^- \ar[d]^{\cdot f} & \sS^s (L^{4e+1} (\alpha_k))
\ar[l]_(0.6){\wrho} \ar[d]^{\Sigma_1} \ar@{-{>>}}[r]^{\eta} & \widetilde \sN (L^{4e+1} (\alpha_k)) \\
\QQ \RhG^+ & \sS^s (L^{4e+3} (\alpha_k)) \ar[l]_(0.6){\wrho}
\ar@{-{>>}}[r]^{\eta} & \widetilde \sN (L^{4e+3} (\alpha_k))
\ar[u]_{\res}^{\cong} }
\end{split}
\end{equation}
and use a similar reasoning.
\end{proof}

\begin{cor}
The function $\wrho \co \sS^s(L^{2d-1}(\alpha_k)) \lra \QQ
\RhG^{(-1)^d}$ is a homomorphism.
\end{cor}

\begin{proof}
Just as in the case $\alpha_1$ it is enough to find in each normal
cobordism class an element such that the addition works for these
representatives. Consider the diagrams (\ref{join-1-k}) and
(\ref{join-k-k}). In the case $d=2e$ we can choose in each normal
cobordism class an element coming from $\sS^s (L^{4e-3} (\alpha_1))$
and in the case $d=2e+1$ choose in each normal cobordism class an
element coming from $\sS^s (L^{4e+1} (\alpha_k))$. Then use the fact
that we have already proved the proposition in the case $d = 2e+2$
and Lemma \ref{f-inverse} just as in the case $\alpha_1$.
\end{proof}

\begin{cor}
The function $\wrho \co \sS^s(L^{2d-1}(\alpha)) \lra \QQ
\RhG^{(-1)^d}$ is a homomorphism.
\end{cor}

The proof is the same as the proof of \cite[Corollary
4.16]{Macko-Wegner(2008)}.

%%%%%%%%%%%%%%%%%%%%%%%%%%%%%%%%%%%%%%%%%

\subsection{Formulas}

\

%%%%%%%%%%%%%%%%%%%%%%%%%%%%%%%%%%%%%%%%%

\

In addition to the information we obtained in the previous
subsection we need formulas to calculate the homomorphism $[\wrho]$
in some cases. The formulas and the proofs are generalizations of
the similar formulas in \cite[section 4.2]{Macko-Wegner(2008)}. We
note that we will need these formulas only in the case $N = 2^K
\cdot M$ for $K \geq 1$ and hence we can assume $2 \; | \; N$ in
this subsection. The starting point is the following:

\begin{thm}{\cite[Theorem 14C.4]{Wall(1999)}} \label{rho-formula-cp}
Let $a = [h \co Q \ra \CC P^{d-1}]$ be an element in $\sS^s (\CC
P^{d-1})$. Then for $t \in S^1$
\[
\wrho_{S^1} (t,a) = \sum_{1 \leq i \leq \lfloor d/2 \rfloor -1} 8
\cdot \bs_{4i} (\eta (a)) \cdot (f^{d-2i} - f^{d-2i-2}) \in \CC,
\]
where $f = (1+t)/(1-t)$.
\end{thm}

This was used to obtain a formula for $\wrho$ for fake lens spaces
which fiber over the fake complex projective spaces in \cite[Theorem
14E.8]{Wall(1999)}. From that we obtain a formula for $[\wrho]$ when
$d=2e$ since in that case there exists in each normal cobordism
class of fake lens spaces a fake lens space which fibers over a fake
complex projective space. When $d=2e+1$ we need to adapt the trick
from \cite[Lemma 4.10, 4.11]{Macko-Wegner(2008)} to our situation.
This is done as follows.

\begin{lem} \label{all-summands-d=2e+1}
Let $d = 2e+1$ and let $c \in \sS^s (L^{2d-1}(\alpha_1))$ be such
that $c = b + a$ where $b = p^{!} (\widetilde b)$ for some
$\widetilde b \in \sS (\CC P^{d-1})$ with $s_{2i} = \bs_{2i} (\eta
\widetilde b)$ and $a$ is such that there exists an $\widetilde a
\in \sN (\CC P^{d-1})$ such that $\bs_{2i} (\widetilde a) = 0$ for
$i < 2e$,  $s_{4e} = \bs_{4e} (\widetilde a)$ and $\eta (a) = p^{!}
(\widetilde a)$. Then
\[
\widetilde\rho (c) = \sum_{1 \leq i \leq \lfloor d/2 \rfloor -1}
\!\! 8 \cdot s_{4i} \cdot (f^{d-2i} - f^{d-2i-2}) + 8 \cdot s_{4e}
\cdot f + z \quad \in \quad \QQ\RhGm
\]
for some $z \in 4 \cdot \RhGm$.
\end{lem}

\begin{proof}
We will show that, if $s_{4i} = 0$ for all $i<e$ and $s_{4e} = 1$
then
\[
\widetilde\rho (a) = 8 f + z \quad \in \quad \QQ\RhGm
\]
for some $z \in 4 \cdot \RhGm$. The proof of the general case is the
same, but the formulas are more complicated.

We will use the suspension map $\Sigma_1$. Our assumptions mean that
$\eta(a)$ is not in the image of the composition $\sS (\CC P^{d-1})
\ra \sN (\CC P^{d-1}) \ra \sN (L^{2d-1}(\alpha_1))$. However,
diagram (\ref{susp-diagram}) and the formula from Theorem
\ref{rho-formula-cp} tells us that $\eta (\Sigma_1(a))$ is in the
image of $\sS (\CC P^{(d+1)-1}) \ra \sN (\CC P^{(d+1)-1}) \ra \sN
(L^{2(d+1)-1}(\alpha_1))$ and we have
\[
f \cdot \widetilde\rho (a) + y = 8 \cdot 1 \cdot (f^2 -1) \quad \in
\quad \QQ\RhGp
\]
for some $y \in 4 \cdot \RhGp$. We obtain the desired identity by
the following calculation. Let $\hrho \in \QQ[\chi]$ and $\hy \in 4
\cdot \ZZ[\chi]$ be representatives for $\wrho (a)$ and $y$. Then
\begin{align*}
(1+\chi)(1-\chi)\hrho + (1-\chi)^2 \hy & \equiv 8 \cdot (4 \chi) \quad \quad \quad \mod \langle 1+ \chi + \cdots + \chi^{N-1} \rangle \\
(1+\chi)(1-\chi)\hrho + (1-\chi)^2 (\hy+8) & \equiv 8 \cdot (1+\chi)^2 \quad \; \mod \langle 1+ \chi + \cdots + \chi^{N-1} \rangle \\
(1+\chi)(1-\chi)\hrho + (1-\chi)^2 (\hy+8) & = 8 \cdot (1+\chi)^2
\!+\! g(\chi) (1 + \chi + \cdots + \chi^{N-1}) \in \QQ[\chi]
\end{align*}
for some $g(\chi) \in \QQ[\chi]$. Because $2$ divides $N$ we get
$\hy+8 = (1+\chi) \cdot w(\chi)$ for some $w(\chi) \in \QQ[\chi]$.
Since $\hy + 8 \in 4 \cdot \ZZ [\chi]$, we obtain $w(\chi) \in 4
\cdot \ZZ[\chi]$. Further write $g (\chi) = 2r+(1+\chi) g'(\chi) = r
(1-\chi) + (1+\chi) (r+g'(\chi))$ for $r \in \QQ$, $g' (\chi) \in
\QQ[\chi]$. We have
\begin{align*}
(1-\chi) \cdot \hrho + (1-\chi)^2 \cdot w(\chi) & = 8 \cdot (1+\chi)
+ g(\chi) \cdot (1 + \chi^2 + \cdots \chi^{N-2}) \in \QQ[\chi]
\end{align*}
and further modulo $\langle 1+ \chi + \cdots + \chi^{N-1} \rangle$
\begin{align*}
(1-\chi) \cdot \hrho + (1-\chi)^2 \cdot w(\chi) & \equiv 8 \cdot (1+\chi) + r \cdot (1-\chi) \cdot (1 + \chi^2 + \cdots \chi^{N-2}) \\
\hrho + (1-\chi) \cdot w(\chi) & \equiv 8 \cdot f + r \cdot (1 +
\chi^2 + \cdots \chi^{N-2})
\end{align*}
Now $(1-\chi) \cdot w(\chi) = (2-(1+\chi)) \cdot w(\chi) = 2 \cdot
w(\chi) - (\hy + 8)$. Further $2 \cdot w(\chi) = w^+(\chi) +
w^-(\chi)$, where $w^\pm(\chi) \in 4 \cdot \RhG^{(\pm 1)}$. Hence
\begin{equation*}
\wrho (a) - 8 \cdot f + w^-(\chi) = (\hy+8) - w^+ (\chi) + r \cdot
(1 + \chi^2 + \cdots + \chi^{N-2})
\end{equation*}
in $\QQ[\chi] / \langle 1+ \chi + \cdots + \chi^{N-1} \rangle$,
while the left hand side of the equation lies in the
$(-1)$-eigenspace and the right-hand side lies in the
$(+1)$-eigenspace and hence both are equal to $0$. It follows that
\[
\wrho (a)= 8 \cdot f - w^-(\chi).
\]
Putting $z = - w^-(\chi)$ yields is the desired formula.
\end{proof}

\begin{prop} \label{rho-formula-lens-sp}
Let $t = (t_{2i})_i \in \widetilde \sN (L^{2d-1}(\alpha_1))_{(2)}$
and $\bar t_{4i} \in \ZZ$ $(1 \leq i \leq c)$ with $\bar t_{4i}
\equiv t_{4i} \; \mod \, 2^K$ and $\bar t_{4i} \equiv 0 \; \mod \,
M^c$. Then we have for the homomorphism $[\wrho] \co \widetilde \sN
(L^{2d-1}(\alpha_1))_{(2)} \lra \QQ \RhG^{(-1)^d} / 4 \cdot
\RhG^{(-1)^d}$ that
\begin{align*}
d = 2e \; : \; [\widetilde \rho] (t) & = \sum_{i=1}^{e-1} 8 \cdot
\bar t_{4i} \cdot f^{d-2i-2} \cdot (f^2-1) \\
d = 2e+1 \; : \; [\widetilde \rho] (t) & = \sum_{i=1}^{e-1} 8 \cdot
\bar t_{4i} \cdot f^{d-2i-2} \cdot (f^2-1) + 8 \cdot \bar t_{4e}
\cdot f.
\end{align*}
\end{prop}

\begin{proof}
It is enough to find for each $t \in \widetilde \sN (L^{2d-1}
(\alpha_1))_{(2)}$ some $a \in  \sS^s (L^{2d-1}(\alpha_1))$ with
$\bt(\eta(a))=(t,0)$ and for which we can calculate $\wrho (a)$. If
$d = 2e$ then there is for each normal cobordism class a fake lens
space which fibers over a fake complex projective space and Theorem
\ref{rho-formula-cp} gives the desired formula. If $d = 2e+1$ we
apply instead Lemma \ref{all-summands-d=2e+1} to get the desired
formula.
\end{proof}

\begin{prop} \label{rho-formula-lens-sp-k}
Let $t = (t_{2i})_i \in \widetilde \sN (L^{2d-1}(\alpha_k))_{(2)}$
and $\bar t_{4i} \in \ZZ$ $(1 \leq i \leq c)$ with $\bar t_{4i}
\equiv t_{4i} \; \mod \, 2^K$ and $\bar t_{4i} \equiv 0 \; \mod \,
M^c$. Then we have for the homomorphism $[\wrho] \co \widetilde \sN
(L^{2d-1}(\alpha_k))_{(2)} \lra \QQ \RhG^{(-1)^d} / 4 \cdot
\RhG^{(-1)^d}$ that
\begin{align*}
d = 2e \; : \; [\widetilde \rho] (t) & = \sum_{i=1}^{e-1} 8 \cdot
\bar t_{4i} \cdot f'_k \cdot f^{d-2i-2} \cdot (f^2-1) \\
d = 2e+1 \; : \; [\widetilde \rho] (t) & = \sum_{i=1}^{e-1} 8 \cdot
\bar t_{4i} \cdot f'_k \cdot f^{d-2i-2} \cdot (f^2-1) + 8 \cdot \bar
t_{4e} \cdot f'_k \cdot f.
\end{align*}
\end{prop}

\begin{proof}
We will use the calculation for $\alpha_1$ and the homeomorphisms
\[
L^{2d+1}(\alpha_k) \cong L^{2d-1}(\alpha_1) \ast L^1(\alpha_k) \quad
\mathrm{and} \quad L^{2d+1}(\alpha_k) \cong L^{2d-1}(\alpha_k) \ast
L^1(\alpha_1)
\]
For $d=2e$ recall the diagram
\begin{equation}
\begin{split}
\xymatrix{ \QQ \RhG^- \ar[d]^{\cdot f_k} & \sS^s (L^{4e-3}
(\alpha_1)) \ar[l]_(0.6){\wrho} \ar[d]^{\Sigma_k} \ar@{-{>>}}[r]^{\eta} & \widetilde \sN (L^{4e-3} (\alpha_1))_{(2)} \\
\QQ \RhG^+ & \sS^s (L^{4e-1} (\alpha_k)) \ar[l]_(0.6){\wrho}
\ar@{-{>>}}[r]^{\eta} & \widetilde \sN (L^{4e-1} (\alpha_k))_{(2)}
\ar[u]_{\res}^{\cong} }
\end{split}
\end{equation}
Let $t \in \sN (L^{4e-1}(\alpha_k))$. Choose $x \in
\sS^s(L^{4e-3}(\alpha_1))$ such that $\bt (\eta (x)) = t = \res
(t)$. Then we have $\bt (\eta (\Sigma_k (x))) = t$ and $[\wrho]
(\eta (\Sigma_k (x))) = [\wrho (x) \cdot f_k]$ can be calculated
using the formulas from the case $k=1$.

For $d=2e+1$ recall the diagram
\begin{equation}
\begin{split}
\xymatrix{ \QQ \RhG^- \ar[d]^{\cdot f} & \sS^s (L^{4e+1} (\alpha_k))
\ar[l]_(0.6){\wrho} \ar[d]^{\Sigma_1} \ar@{-{>>}}[r]^{\eta} & \widetilde \sN (L^{4e+1} (\alpha_k))_{(2)} \\
\QQ \RhG^+ & \sS^s (L^{4e+3} (\alpha_k)) \ar[l]_(0.6){\wrho}
\ar@{-{>>}}[r]^{\eta} & \widetilde \sN (L^{4e+3} (\alpha_k))_{(2)}
\ar[u]_{\res}^{\cong} }
\end{split}
\end{equation}
Let $t \in \widetilde \sN (L^{4e+1}(\alpha_k))$. Choose $x \in
\sS^s(L^{4e+1}(\alpha_1))$ such that $\bt (\eta (x)) = t$. Then we
have $\bt (\eta (\Sigma_1 (x))) = t$ and $\wrho (\Sigma_1 (x)) =
\wrho (x) \cdot f$. We obtain the equation
\[
f \cdot \wrho(x) + y = \sum_{i=1}^{e-1} 8 \cdot \bar t_{4i} \cdot
f'_k \cdot f^{d+1-2i-2} \cdot (f^2-1) + 8 \cdot \bar t_{4e} \cdot
f'_k \cdot (f^2 -1) \in  \QQ \RhG^+
\]
for some $y \in 4 \cdot \RhG^+$ using the formulas from the case $d
= 2e+2$ which we have already dealt with. Now a modification of the
argument from the proof of Lemma \ref{all-summands-d=2e+1} can be
used to obtain the formula for $[\wrho] (\eta (x))$.
\end{proof}

%%%%%%%%%%%%%%%%%%%%%%%%%%%%%%%%%%%%%%%%%%%%

\section{Revision} \label{sec:revision}

%%%%%%%%%%%%%%%%%%%%%%%%%%%%%%%%%%%%%%%%%%%%

In this section we show how to use Proposition \ref{ses-vs-rep-thy}
to calculate the structure set of the fake lens spaces. We also
remind the reader the calculations when $N = M$ odd and $N = 2^K$
which are used below to solve the general case.

\begin{prop} \label{prop:sset-general}
Let $d \geq 3$ and $N = 2^K \cdot M$ with $K \geq 1$, $M$ odd, and $G = \ZZ_N$. Then we have
\[
 \sS^s (L_N^{2d-1}(\alpha)) \cong \bar \Sigma_N(d) \oplus \ker [\wrho_N]
\]
with $\bar \Sigma_N (d) := \im (\wrho_N \co \sS^s (L(\alpha)) \ra
\QQ \RhG^{(-1)^d})$. The rank of the free abelian group $\bar
\Sigma_N (d)$ equals the rank of $\wL^s_{2d} (G)$.
\end{prop}

The proof is the same as the proof of Theorem 5.1 in
\cite{Macko-Wegner(2008)}. So in order to calculate the structure
set one needs to calculate the kernel $\ker [\wrho_N]$. This is
essentially what Wall has done in \cite[chapter 14E]{Wall(1999)} and
what we have done in \cite{Macko-Wegner(2008)} in the special cases.

\begin{thm} \textup{(}\cite[chapter 14E]{Wall(1999)}\textup{)}\label{Wall-thm}
Let $d \geq 3$. If $N=M$ is odd then we have
\[
\sS^s (L_M^{2d-1}(\alpha)) \cong \bar \Sigma_M (d).
\]
Equivalently, $\ker [\wrho_M] = 0$.
\end{thm}

\begin{thm} \textup{(}\cite[Theorem 1.2]{Macko-Wegner(2008)}\textup{)}\label{main-thm-2-K}
Let $d \geq 3$. If $N=2^K$ with $K \geq 0$ then we have
\[
\sS^s (L_{2^K}^{2d-1}(\alpha)) \cong \bar \Sigma_{2^K}(d) \oplus \ker
[\wrho_{2^K}] \cong \bar \Sigma_{2^K}(d) \oplus \bigoplus_{i=1}^{c}
\ZZ_{2^{\min\{K,1\}}} \oplus \bigoplus_{i=1}^{c}
\ZZ_{2^{\min\{K,2i\}}}
\]
where $\bar \Sigma_{2^K}(d)$ is a free abelian group of rank $N/2-1$
if $d=2e+1$ and $N/2$ if $d=2e$ and $c = \lfloor (d-1)/2 \rfloor$.
\end{thm}

%%%%%%%%%%%%%%%%%%%%%%%%%%%%%%%%%%%%%%%%%%%%

\section{Calculations} \label{sec:calculations}

%%%%%%%%%%%%%%%%%%%%%%%%%%%%%%%%%%%%%%%%%%%%

In this section we give the proof of Theorem \ref{main-thm} by
calculating $\ker [\wrho_N]$. Notice that for any fake lens space
$L^{2d-1}_N(\alpha)$ there exists $k \in \NN$ and a homotopy
equivalence $h \co L^{2d-1}_N(\alpha) \lra L^{2d-1}_N(\alpha_k)$
(see for instance \cite[Corollary 2.4]{Macko-Wegner(2008)}). It
induces an isomorphism $h_* \co \sS^s (L^{2d-1}_N(\alpha)) \to \sS^s
(L^{2d-1}_N(\alpha_k))$. Hence it suffices to consider the case
$\alpha = \alpha_k$. In Proposition \ref{ker-prop} we calculate
\[
 \ker [\wrho_N]: \widetilde \sN (L^{2d-1}_N(\alpha_k)) \to \QQ R_{\widehat \ZZ_N}^{(-1)^d}/4 \cdot R_{\widehat \ZZ_N}^{(-1)^d}.
\]
Theorem \ref{main-thm} follows from this proposition together with
Proposition \ref{prop:sset-general} and Theorems \ref{Wall-thm},
\ref{main-thm-2-K}.

We need a little preparation for Proposition \ref{ker-prop}. Notice
that we have commutative diagrams
\begin{equation}\label{diagram_2^K}
\begin{split}
\xymatrix{
\widetilde \sN (L^{2d-1}_N(\alpha)) \ar[d]_{(p_{2^K}^N)^{!}} \ar[r]^(0.435){[\wrho_N]} & \QQ R_{\widehat \ZZ_N}^{(-1)^d}/4 \cdot R_{\widehat \ZZ_N}^{(-1)^d} \ar[d] \\
\widetilde \sN (L^{2d-1}_{2^K}(\alpha)) \ar[r]^(0.435){[\wrho_{2^K}]} & \QQ
R_{\widehat \ZZ_{2^K}}^{(-1)^d}/4 \cdot R_{\widehat
\ZZ_{2^K}}^{(-1)^d} }
\end{split}
\end{equation}
and
\begin{equation}\label{diagram_M}
\begin{split}
\xymatrix{
\widetilde \sN (L^{2d-1}_N(\alpha)) \ar[d]_{(p_M^N)^{!}} \ar[r]^(0.435){[\wrho_N]} & \QQ R_{\widehat \ZZ_N}^{(-1)^d}/4 \cdot R_{\widehat \ZZ_N}^{(-1)^d} \ar[d] \\
\widetilde \sN (L^{2d-1}_M(\alpha)) \ar[r]^(0.435){[\wrho_M]} & \QQ
R_{\widehat \ZZ_M}^{(-1)^d}/4 \cdot R_{\widehat \ZZ_M}^{(-1)^d}. }
\end{split}
\end{equation}
This follows from the fact that the $\wrho$-invariant is natural
with respect to the restriction maps. Note also that if we use our
identification of the right hand terms with the polynomial rings,
then the right hand vertical maps are given by $\chi \mapsto \chi$,
see Remark \ref{chi-rem}. We also have an isomorphism
\[
\widetilde \sN (L_N ) \cong \widetilde \sN (L_{N} )_{(2)} \oplus
\widetilde \sN (L_N )_{(odd)} \xra{(p_{2^K}^N)^{!} \oplus
(p_M^N)^{!}} \widetilde \sN (L_{2^K} ) \oplus \widetilde \sN (L_M )
\]
where $L_? = L^{2d-1}_? (\alpha)$, see (\ref{ni-N-vs-2K-M}).

\begin{prop} \label{ker-prop}
Let $N = 2^K \cdot M$ with $K \geq 1$, $M>1$ odd and let $d \geq 3$.
Then we have
\begin{align*}
\ker [\wrho_N] & \cong \big( (p_{2^K}^N)^{!} \big)^{-1} \ker [\wrho_{2^K}] \oplus \big( (p_{M}^N)^{!}\big)^{-1} \ker [\wrho_M] = \big( (p_{2^K}^N)^{!} \big)^{-1} \ker [\wrho_{2^K}] \oplus 0 \\
& \subseteq \quad \widetilde \sN (L^{2d-1}_N(\alpha_k))_{(2)} \oplus
\widetilde \sN (L^{2d-1}_N(\alpha_k))_{(odd)}.
\end{align*}
\end{prop}

For the proof of this proposition we need the following two lemmas.

\begin{lem}\label{decomposition-lem}
Let $a \in \QQ[\chi]$ such that there exist $b, c \in 4 \cdot
\ZZ[\chi]$ satisfying
\begin{enumerate}
 \item $a \equiv b \quad \mod 1 + \chi + \cdots + \chi^{2^K-1}$ \quad and
 \item $a \equiv c \quad \mod 1 + \chi^{2^K} + \cdots + \chi^{2^K \cdot (M-1)}$.
\end{enumerate}
Then there exists $d \in 4 \cdot \ZZ[\chi]$ such that
\[
 M \cdot a \equiv d \quad \mod 1 + \chi + \cdots + \chi^{N-1}
\]
holds.
\end{lem}

\begin{proof}
Our candidate for $d$ is
\[
 d := c \cdot M + (b - c) \cdot \left( 1 + \chi^{2^K} + \cdots + \chi^{2^K \cdot (M-1)} \right).
\]
We have to check that $M \cdot a \equiv d \; \mod 1 + \chi + \cdots
+ \chi^{N-1}$ holds. There is the factorization
\[
 1 + \chi + \cdots + \chi^{N-1} = (1 + \chi + \cdots + \chi^{2^K-1}) \cdot (1 + \chi^{2^K} + \cdots + \chi^{2^K \cdot (M-1)}).
\]
Notice that the two factors are prime to each other in $\QQ[\chi]$
because of
\begin{equation}\label{M-eq}
 \sum_{i=0}^{M-1} \left( \chi^{2^K} \right)^i \equiv \sum_{i=0}^{M-1} 1^i \equiv M \quad \mod 1 + \chi + \cdots + \chi^{2^K-1}.
\end{equation}
Therefore, it suffices to check that
\begin{eqnarray*}
M \cdot a & \equiv & d \quad \mod 1 + \chi + \cdots + \chi^{2^K-1} \qquad \mbox{and} \\
M \cdot a & \equiv & d \quad \mod 1 + \chi^{2^K} + \cdots +
\chi^{2^K \cdot (M-1)}
\end{eqnarray*}
hold. But this is true because of equation \ref{M-eq}.
\end{proof}

\begin{lem}\label{M-factor-lem}
Let $k \in \NN$ with $(k,N)=1$ and $q \in \ZZ[x]$. Then
\begin{eqnarray*}
 8 \cdot f'_k \cdot (f^2-1) \cdot M^{2+2\deg(q)} \cdot q(f^2) & \in & 4 \cdot \ZZ[\chi] / \langle 1 + \chi^{2^K} + \cdots + \chi^{2^K \cdot (M-1)} \rangle \; \mbox{and}\\
 8 \cdot f'_k \cdot f \cdot M^{1+2\deg(q)} \cdot q(f^2) & \in & 4 \cdot \ZZ[\chi] / \langle 1 + \chi^{2^K} + \cdots + \chi^{2^K \cdot (M-1)} \rangle.
\end{eqnarray*}
\end{lem}

\begin{proof}
Because of the equation
\begin{align*}
 & (1-\chi) \cdot \left( 1 + \chi + \cdots + \chi^{2^K-1} \right) \cdot \left( 1 + 2 \cdot \chi^{2^K} + 3 \cdot \chi^{2 \cdot 2^K} + \cdots + M \cdot \chi^{(M-1) \cdot 2^K} \right) =\\
 & \left( 1 + \chi^{2^K} + \cdots + \chi^{(M-1) \cdot 2^K} \right) - M \cdot \chi^N
\end{align*}
we have in $\QQ[\chi] / \langle 1 + \chi^{2^K} + \cdots + \chi^{2^K
\cdot (M-1)} \rangle$:
\[
 M \cdot f = - (1+\chi) \cdot \left( 1 + \chi + \cdots + \chi^{2^K-1} \right) \cdot \left( 1 + 2 \chi^{2^K} + 3 \chi^{2 \cdot 2^K} + \cdots + M \chi^{(M-1) \cdot 2^K} \right).
\]
Notice that $f'_k \in \ZZ[\chi] / \langle 1 + \chi^{2^K} + \cdots +
\chi^{2^K \cdot (M-1)} \rangle$ because of $f'_k \in R_{\widehat
\ZZ_N}$. Therefore, we obtain
\begin{eqnarray*}
 8 \cdot f'_k \cdot (f^2-1) \cdot M^{2+2\deg(q)} \cdot q(f^2) & \in & 4 \cdot \ZZ[\chi] / \langle 1 + \chi^{2^K} + \cdots + \chi^{2^K \cdot (M-1)} \rangle \; \mbox{and}\\
 8 \cdot f'_k \cdot f \cdot M^{1+2\deg(q)} \cdot q(f^2) & \in & 4 \cdot \ZZ[\chi] / \langle 1 + \chi^{2^K} + \cdots + \chi^{2^K \cdot (M-1)} \rangle.
\end{eqnarray*}
\end{proof}

\begin{proof}[Proof of Proposition \ref{ker-prop}]
In the sequel we use the identification
\[
\widetilde \sN (L^{2d-1}_N(\alpha_k)) = \widetilde \sN
(L^{2d-1}_N(\alpha_k))_{(2)} \oplus \widetilde \sN
(L^{2d-1}_N(\alpha_k))_{(odd)}.
\]
Obviously, we have
\[
 \ker [\wrho_N] \supseteq \ker([\wrho_N]|_{\widetilde \sN (L^{2d-1}_{N} (\alpha_k))_{(2)}}) \oplus \ker([\wrho_N]|_{\widetilde \sN (L^{2d-1}_N (\alpha_k))_{(odd)}}).
\]
We want to prove the other inclusion as well. Let $x = x_{(2)} +
x_{(odd)} \in \ker [\wrho_N]$ with $x_{(2)} \in \widetilde \sN
(L^{2d-1}_{N} (\alpha_k))_{(2)}$ and $x_{(odd)} \in \widetilde \sN
(L^{2d-1}_N (\alpha_k))_{(odd)}$. Notice that $2^K \cdot x_{(2)} =
0$ and $M^c \cdot x_{(odd)} = 0$ hold. Since $M$ is odd, there exist
$a,b \in \ZZ$ such that $a \cdot 2^K + b \cdot M^c = 1$. We conclude
\[
 [\wrho_N](x_{(2)}) = [\wrho_N](x_{(2)} - a \cdot 2^K \cdot x_{(2)}) = [\wrho_N](b \cdot M^c \cdot x_{(2)}) = [\wrho_N](b \cdot M^c \cdot x) = 0
\]
and
\[
 [\wrho_N](x_{(odd)}) = [\wrho_N](x_{(odd)} - b \cdot M^c \cdot x_{(odd)}) = [\wrho_N](a \cdot 2^K \cdot x_{(odd)}) = [\wrho_N](a \cdot 2^K \cdot x) = 0.
\]
This shows
\[
 \ker [\wrho_N] = \ker([\wrho_N]|_{\widetilde \sN (L^{2d-1}_{N} (\alpha_k))_{(2)}}) \oplus \ker([\wrho_N]|_{\widetilde \sN (L^{2d-1}_N (\alpha_k))_{(odd)}}).
\]
It remains to prove that we have $\ker([\wrho_N]|_{\widetilde \sN
(L^{2d-1}_{N} (\alpha_k))_{(2)}}) = \big( (p_{2^K}^N)^{!} \big)^{-1}
\ker [\wrho_{2^K}]$ and $\ker([\wrho_N]|_{\widetilde \sN (L^{2d-1}_N
(\alpha_k))_{(odd)}}) = \big( (p_{M}^N)^{!}\big)^{-1} \ker [\wrho_M]
= 0$. From diagram \ref{diagram_M} we conclude
$\ker([\wrho_N]|_{\widetilde \sN (L^{2d-1}_N (\alpha_k))_{(odd)}})
\subseteq \big( (p_{M}^N)^{!}\big)^{-1} \ker [\wrho_M]$. Proposition
\ref{ses-vs-rep-thy} implies that $\eta \co
\sS^s(L^{2d-1}_M(\alpha_k)) \to \widetilde
\sN(L^{2d-1}_M(\alpha_k))$ induces an isomorphism $\ker(\wrho_M)
\cong \ker [\wrho_M]$. Notice that $\ker [\wrho_M] = 0$ by Theorem
\ref{Wall-thm}. Hence we obtain
\[
 \ker([\wrho_N]|_{\widetilde \sN (L^{2d-1}_N (\alpha_k))_{(odd)}})) = \big( (p_{M}^N)^{!}\big)^{-1} \ker [\wrho_M] = 0.
\]
Diagram \ref{diagram_2^K} implies
\[
 \ker([\wrho_N]|_{\widetilde \sN (L^{2d-1}_{N} (\alpha_k))_{(2)}}) \subseteq \big( (p_{2^K}^N)^{!} \big)^{-1} \ker [\wrho_{2^K}].
\]
The remaining part of the proof of $\ker([\wrho_N]|_{\widetilde \sN
(L^{2d-1}_{N} (\alpha_k))_{(2)}}) \supseteq \big( (p_{2^K}^N)^{!}
\big)^{-1} \ker [\wrho_{2^K}]$ needs more effort. Let $x_{(2)} \in
\big( (p_{2^K}^N)^{!} \big)^{-1} \ker [\wrho_{2^K}] \subseteq
\widetilde \sN (L^{2d-1}_{N}(\alpha_k))_{(2)}$. Recall from
(\ref{ni-N-vs-2K-M}) that we have isomorphisms
\begin{align*}
 \bt_{(2)} = (\bt_{2i})_i & \co \widetilde \sN (L^{2d-1}_{N}(\alpha_k))_{(2)} \xra{\cong} \bigoplus_{i=1}^c \ZZ_{2^K} \oplus \bigoplus_{i=1}^c \ZZ_2 \\
 \bt = (\bt_{2i})_i & \co \widetilde \sN (L^{2d-1}_{2^K}(\alpha_k)) \xra{\cong} \bigoplus_{i=1}^c \ZZ_{2^K} \oplus \bigoplus_{i=1}^c \ZZ_2
\end{align*}
which commute with the isomorphism $(p_{2^K}^N)^{!}$. Choose $\bar
t_{4i} \in \ZZ$ $(1 \leq i \leq c)$ with $\bar t_{4i} \equiv
\bt_{4i}\big((p_{2^K}^N)^{!}(x_{(2)})\big) \; \mod \, 2^K$. From
Proposition \ref{rho-formula-lens-sp-k} we conclude
\begin{equation*}
 [\wrho_{2^K}]\big((p_{2^K}^N)^{!}(x_{(2)})\big) = \left\{ \begin{array}{ll}
 8 \cdot f'_k \cdot (f^2-1) \cdot q_{\bar t}(f^2) & \quad \mbox{$d$ even},\\
 \rule{0in}{3ex}
 8 \cdot f'_k \cdot f \cdot q_{\bar t}(f^2) & \quad  \mbox{$d$ odd}
 \end{array} \right.
\end{equation*}
where the polynomial $q_{\bar t} \in \ZZ[x]$ is defined by
\begin{equation*}
 q_{\bar t}(x) := \left\{ \begin{array}{ll}
 \displaystyle \sum_{i=0}^{c-1} {\bar t}_{4(i+1)} \cdot x^{c-i-1} & \quad \mbox{$d$ even},\\
 \rule{0in}{3ex}
 \displaystyle \sum_{i=1}^{c-1} {\bar t}_{4i} \cdot x^{c-i-1} \cdot (x-1) + {\bar t}_{4c} & \quad \mbox{$d$ odd}.
 \end{array} \right.
\end{equation*}
Since $x_{(2)} \in \ker [\wrho_{2^K}]$, we have
\begin{eqnarray*}
 8 \cdot f'_k \cdot (f^2-1) \cdot q_{\bar t}(f^2) & \in & 4 \cdot R_{\widehat \ZZ_{2^K}} \cong 4 \cdot \ZZ[\chi] / \langle 1 + \chi + \cdots + \chi^{2^K-1} \rangle \quad \mbox{($d$ even)},\\
 8 \cdot f'_k \cdot f \cdot q_{\bar t}(f^2) & \in & 4 \cdot R_{\widehat \ZZ_{2^K}} \cong 4 \cdot \ZZ[\chi] / \langle 1 + \chi + \cdots + \chi^{2^K-1} \rangle \quad \mbox{($d$ odd)}.
\end{eqnarray*}
Lemma \ref{M-factor-lem} tells us in respective cases that
\begin{eqnarray*}
 8 \cdot f'_k \cdot (f^2-1) \cdot M^{2+2\deg(q_{\bar t})} \cdot q_{\bar t}(f^2) & \in & 4 \cdot \ZZ[\chi] / \langle 1 + \chi^{2^K} + \cdots + \chi^{2^K \cdot (M-1)} \rangle \\
 8 \cdot f'_k \cdot f \cdot M^{1+2\deg(q_{\bar t})} \cdot q_{\bar t}(f^2) & \in & 4 \cdot \ZZ[\chi] / \langle 1 + \chi^{2^K} + \cdots + \chi^{2^K \cdot (M-1)} \rangle.
\end{eqnarray*}
Using Lemma \ref{decomposition-lem} we obtain in respective cases
\begin{eqnarray*}
 8 \cdot f'_k \cdot (f^2-1) \cdot M^{3+2\deg(q_{\bar t})} \cdot q_{\bar t}(f^2) & \in & 4 \cdot \ZZ[\chi] / \langle 1 + \chi + \cdots +\chi^{N-1} \rangle \cong 4 \cdot R_{\widehat \ZZ_N} \\
 8 \cdot f'_k \cdot f \cdot M^{2+2\deg(q_{\bar t})} \cdot q_{\bar t}(f^2) & \in & 4 \cdot \ZZ[\chi] / \langle 1 + \chi + \cdots + \chi^{N-1} \rangle \cong 4 \cdot R_{\widehat \ZZ_N}.
\end{eqnarray*}
Let $z \in \ZZ$ with $z \cdot M^{3+2\deg(q_t)+c} \equiv 1 \; \mod
2^K$. We define $\bar t'_{4i} := z \cdot M^{3+2\deg(q_{\bar t})+c}
\cdot \bar t_{4i} \in \ZZ$ $(1 \leq i \leq c)$ and conclude
\[
 \bar t'_{4i} \equiv \bar t_{4i} \equiv \bt_{4i}\big((p_{2^K}^N)^{!}(x_{(2)})\big) = \bt_{4i}(x_{(2)}) \; \mod \, 2^K, \quad
 \bar t'_{4i} \equiv 0 \; \mod \, M^c.
\]
Proposition \ref{rho-formula-lens-sp-k} tells us that $[\widetilde
\rho_N] (x_{(2)}) \in \QQ R_{\widehat \ZZ_N}^{(-1)^d}/4 \cdot
R_{\widehat \ZZ_N}^{(-1)^d}$ is given by
\begin{align*}
d = 2e \; : \; [\widetilde \rho_N] (x_{(2)}) & = \sum_{i=1}^{e-1} 8
\cdot
\bar t'_{4i} \cdot f'_k \cdot f^{d-2i-2} \cdot (f^2-1) \\
& = z \cdot M^c \cdot 8 \cdot f'_k \cdot (f^2-1) \cdot M^{3+2\deg(q_{\bar t})} \cdot q_{\bar t}(f^2) \in 4 \cdot R_{\widehat \ZZ_N} \\
d = 2e+1 \; : \; [\widetilde \rho_N] (x_{(2)}) & = \sum_{i=1}^{e-1}
8 \cdot
\bar t'_{4i} \cdot f'_k \cdot f^{d-2i-2} \cdot (f^2-1) + 8 \cdot \bar t'_{4e} \cdot f'_k \cdot f \\
& = z \cdot M^{c+1} \cdot 8 \cdot f'_k \cdot f \cdot
M^{2+2\deg(q_{\bar t})} \cdot q_{\bar t}(f^2) \in 4 \cdot
R_{\widehat \ZZ_N}.
\end{align*}
Hence we have $x_{(2)} \in \ker([\wrho_N]|_{\widetilde \sN
(L^{2d-1}_{N} (\alpha_k))_{(2)}})$.
\end{proof}

%%%%%%%%%%%%%%%%%%%%%%%%%%%%%%%%%%%%%%%%%%%%%%%%%

\section{The suspension homomorphism} \label{sec:suspension}

%%%%%%%%%%%%%%%%%%%%%%%%%%%%%%%%%%%%%%%%%%%%%%%%%

The suspension homomorphism
\begin{equation} \label{eq:suspension}
 \Sigma \colon \sS^s (L^{2d-1} (\alpha_k)) \ra \sS^s (L^{2d+1}
 (\alpha_k))
\end{equation}
was already mentioned in section \ref{sec:ses}, see also
\cite[chapter 14]{Wall(1999)}.  It is an interesting question on its
own to understand this map. It was also used as an important tool in \cite[chapter 14]{Wall(1999)} and \cite{LdM(1971)} to obtain the calculation of the structure sets of lens spaces when $N = 2$ and $N$ odd. We were able to perform our calculations essentially without the use of this map. On the other hand, understanding of $\Sigma$ enables us to give a description of the torsion invariants of fake lens spaces in the next section.

Studying $\Sigma$ is a problem closely related to studying the
splitting of simple homotopy equivalences along submanifolds. Recall
that a simple homotopy equivalence $h \co M^n \ra X^n$ of manifolds
is called {\it split} along a locally flat submanifold $Y^{n-q}
\subset X^n$ if it is transverse to it and the restrictions $h| \co
h^{-1} (Y) \ra Y$ and $h| \co M \smallsetminus h^{-1} (Y) \ra X
\smallsetminus Y$ are simple homotopy equivalences. We also say that
$h$ {\it can be made split} if it is homotopic to a split map. The
question whether a simple homotopy equivalence can be made split is
called a {\it splitting problem}. It may be obstructed as we discuss
below.

Before that notice that if $h \co L \ra L^{2d-1} (\alpha_k)$ is a
simple homotopy equivalence representing an element $x \in \sS^s
(L^{2d-1} (\alpha_k))$ then the suspension $\Sigma (x) \in \sS^s
(L^{2d+1} (\alpha_k))$ is split along $L^{2d-1} (\alpha_k)$. The
converse is also true \cite{Wall(1999)}. Therefore studying $\Sigma$
is equivalent to studying the splitting problems with the target
$L^{2d+1} (\alpha_k)$ along the embedded $L^{2d-1} (\alpha_k)
\subset L^{2d+1} (\alpha_k)$.

There is a general obstruction theory for the splitting problems,
the obstruction groups are the so-called $LS$-groups which are
renamed $LN$-groups in the special case when the inclusion $Y \subset
X$ induces an isomorphism on the fundamental groups. They depend on
the homomorphism $\pi_1 (X \smallsetminus Y) \ra \pi_1 (X)$, the dimension $n-q$ and the codimension $q$ (which is usually not indicated in the notation). In fact, in our case $L^{2d-1} (\alpha_k) \subset L^{2d+1} (\alpha_k)$ there is an exact sequence (see \cite[section 7.2]{Ranicki(1981)}):
\[
\cdots \ra LN_{2d} (\ZZ \ra \ZZ_N) \ra \sS^s (L^{2d-1} (\alpha_k))
\xra{\Sigma} \sS^s (L^{2d+1} (\alpha_k)) \ra LN_{2d-1} (\ZZ \ra
\ZZ_N) \ra \cdots
\]
The $LN$-groups were extensively studied in \cite{LdM(1971)},
\cite{Wall(1999)}, \cite{Ranicki(1981)}. However, we will not
directly need these calculations. We will only use the calculations
of \cite{Wall(1999)} and \cite{LdM(1971)} when $N = 2$ and $N$ is odd. These calculations are obtained using the $LN$-groups $LN_\ast (\ZZ \ra
\ZZ_{N'})$ when $N'=2$, so in fact we use them indirectly. The point
is that in order to obtain results for general $N \geq 2$ only the
groups $LN_\ast (\ZZ \ra \ZZ_{N'})$ are needed where $N' = 2$. In
particular we do not need the complicated calculation of
\cite[section 7.8]{Ranicki(1981)}.

In Remark \ref{connection-to-BL} we also describe a relation of
invariants we obtain to the so-called Browder-Livesay invariants.

\begin{thm} \label{thm1}
Let $G = \ZZ_N$ with $N = 2^K \cdot M$, $M$ odd and $K \geq 1$ and
$e \geq 1$. Then we have for the map
\[
 \Sigma \colon \sS^s (L_N^{4e+1} (\alpha_k)) \ra \sS^s (L_N^{4e+3}
 (\alpha_k))
\]
that:
\begin{enumerate}
 \item it is split-injective.
 \item $\sS^s (L_N^{4e+3} (\alpha_k)) \cong \im
 (\Sigma) \oplus \ZZ (\sigma)$ where $\sigma \in
 \widetilde{L}^s_{4e+4} (G) \subset \sS^s (L_N^{4e+3}
 (\alpha_k))$ is defined by $\wrho (\sigma) = 8 \in 4
 \cdot \RhG^+$.
 \item $\im  (\Sigma) = \ker \big( \sS^s (L_N^{4e+3} (\alpha_k))
\xra{\wrho} \QQ \RhG \xra{\chi \mapsto -1} \QQ \big)$.
\end{enumerate}
\end{thm}

\begin{thm} \label{thm2}
Let $G = \ZZ_N$ with $N = 2^K \cdot M$, $M$ odd and $K \geq 1$, $e
\geq 2$. Then we have for the map
\[
 \Sigma \colon \sS^s (L_N^{4e-1} (\alpha_k)) \ra \sS^s (L_N^{4e+1} (\alpha_k))
\]
that:
\begin{enumerate}
 \item $\ker (\Sigma) = \ZZ(\omega)$ where $\omega \in \wL^s_{4e}(\ZZ_N) \subset \sS^s (L_N^{4e-1} (\alpha_k))$ is defined by
 \[
 \wrho (\omega) = 16 \cdot (1 + \chi^2 + \cdots + \chi^{N-2}) \in 4 \cdot \RhG^+.
 \]
 \item $\sS^s (L_N^{4e+1} (\alpha_k)) \cong \im (\Sigma) \oplus \ZZ_2 (\mu_{4e-2})$ where $\mu_{4e-2} \in \sS^s (L_N^{4e+1} (\alpha_k))$ is defined by $\wrho (\mu_{4e-2}) = 0$, $\bt_{(odd)} (\eta (\mu_{4e-2}) ) = 0$, $\bt_{4e-2} (\eta (\mu_{4e-2}) ) = 1$ and $\bt_{(2i)} (\eta (\mu_{4e-2}) ) = 0$ for $i \neq 2e-1$.
 \item $\im (\Sigma) = \ker \big( \sS^s (L_N^{4e+1} (\alpha_k)) \xra{\bt_{4e-2} \circ \eta} \ZZ_2 \big)$.
\end{enumerate}
\end{thm}

\begin{rem}
The two theorems can be summarized in the following exact sequences
(when $e \geq 1$ or $e \geq 2$ respectively):
\begin{align*}
0 & \ra \sS^s (L_N^{4e+1} (\alpha_k)) \xra{\Sigma} \sS^s (L_N^{4e+3} (\alpha_k)) \ra \ZZ(\sigma) \ra 0 \\
0 \ra \ZZ(\omega) & \ra \sS^s (L_N^{4e-1} (\alpha_k)) \xra{\Sigma}
\sS^s (L_N^{4e+1} (\alpha_k)) \ra \ZZ_2 (\mu_{4e-2}) \ra 0.
\end{align*}
\end{rem}

\subsection{Proof of Theorem \ref{thm1}}

\begin{lem} \label{restriction-3-to-1}
 The restriction map
 \[
  \res \co \wsN (L_N^{4e+3} (\alpha_k)) \ra \wsN (L_N^{4e+1} (\alpha_k))
 \]
 is an isomorphism.
\end{lem}
\begin{proof}
The following diagram implies that the restriction map is
surjective.
\[
\xymatrix{
 \sS^s (L_N^{4e+1}(\alpha_k)) \ar[d]_{\Sigma} \ar@{-{>>}}[r]^{\eta} & \wsN (L_N^{4e+1}(\alpha_k)) \\
 \sS^s (L_N^{4e+3}(\alpha_k)) \ar@{-{>>}}[r]^{\eta} & \wsN (L_N^{4e+3}(\alpha_k)) \ar[u]_{\res}
}
\]
The proof of injectivity is more difficult. We first consider the
case $N$ odd. By an Atiyah-Hirzebruch spectral sequence argument
Wall analyzed that the order of the group $\wsN (L_N^{2d-1}
(\alpha_k)) \cong KO (L^{2d-1}(\alpha)) \otimes \ZZ[\frac{1}{2}]$ is
$N^c$ with $c = \lfloor (d-1)/2 \rfloor$. Hence we have $\big| \wsN
(L_N^{4e+3} (\alpha_k)) \big| = \big| \wsN (L_N^{4e+1} (\alpha_k))
\big|$ and conclude that the restriction map is bijective for $N$
odd. For $N$ even we have the following diagram.
\[
\xymatrix{
 \sN (L_N^{4e+3}(\alpha_k))\ar[d]_{\res} \quad \cong & \hspace{-0.5cm} \bigoplus_{i=1}^e H^{4i} (L^{4e+3}(\alpha_k);\ZZ) \oplus \bigoplus_{i=1}^e H^{4i-2}(L^{4e+3}(\alpha_k);\ZZ_2) \ar[d] \\
 \sN (L_N^{4e+1}(\alpha_k)) \quad \cong & \hspace{-0.5cm} \bigoplus_{i=1}^e H^{4i} (L^{4e+1}(\alpha_k);\ZZ) \oplus \bigoplus_{i=1}^e H^{4i-2}(L^{4e+1}(\alpha_k);\ZZ_2)
}
\]
The map on the right hand side is induced by the inclusion
$L^{4e+1}(\alpha_k) \subset L^{4e+3}(\alpha_k)$ and is hence an
isomorphism. This shows that $\res \co \wsN (L_N^{4e+3} (\alpha_k))
\ra \wsN (L_N^{4e+1} (\alpha_k))$ is an isomorphism for $N$ even. It
remains to consider the case $N = M \cdot 2^K$ with $M > 1$ odd and
$K \geq 1$. The diagram
\[
\xymatrix{
 \sN (L_N^{4e+3}(\alpha_k))\ar[d]_{\res} \ar[rr]^(0.4){(p^N_M)^! \oplus (p^N_{2^K})^!}_(0.4){\cong} & & \sN (L_M^{4e+3}(\alpha_k)) \oplus \sN (L_{2^K}^{4e+3}(\alpha_k)) \ar[d]_{\res \oplus \res}^\cong\\
 \sN (L_N^{4e+1}(\alpha_k)) \ar[rr]^(0.4){(p^N_M)^! \oplus (p^N_{2^K})^!}_(0.4){\cong} & & \sN (L_M^{4e+1}(\alpha_k)) \oplus \sN (L_{2^K}^{4e+1}(\alpha_k))
}
\]
implies that the restriction map is an isomorphism.
\end{proof}

\begin{proof}[Proof of Theorem 6.1]

\

(1) Let $x \in \ker (\Sigma)$. Then
\[
 \eta (x) = \res \circ \eta \circ \Sigma (x) = 0 \quad \textup{and hence} \quad x \in \wL^s_{4e+2} (G) \subset \sS^s (L_N^{4e+1} (\alpha_k)).
\]
Further
\[
 0 = \wrho (\Sigma (x)) = f \cdot \wrho (x).
\]
By \cite[Lemma 5.6]{Macko-Wegner(2008)} we obtain $\wrho (x) = 0$.
These two facts together imply $x =0$.

(2) First we show
\[
 \sS^s (L_N^{4e+3} (\alpha_k)) \cong \im (\Sigma) + \ZZ (\sigma).
\]
Let $x \in \sS^s (L_N^{4e+3} (\alpha_k))$. Choose $y \in \sS^s
(L_N^{4e+1} (\alpha_k))$ such that $\eta (y) = \res \circ \eta (x)$.
Then we have $\res \circ \eta (x - \Sigma (y)) = 0$ and by Lemma
\ref{restriction-3-to-1} also $\eta (x - \Sigma (y)) = 0$ and hence
$x - \Sigma (y) \in \wL^s_{4e+4} (G) \subset \sS^s (L_N^{4e+3}
(\alpha_k))$.

Let $8 \cdot z$ where $z \in \ZZ$ be the image of $x - \Sigma (y)$
under the composition
\[
 \sS^s (L_N^{4e+3} (\alpha_k)) \xra{\wrho} 4 \cdot \RhG^+ \xra{\chi \mapsto -1} 8 \cdot \ZZ.
\]
Then we have that
\[
 \wrho (x - \Sigma (y) - z \cdot \sigma) \in 4 \cdot \RhG^+ \subset 4 \cdot \ZZ[\chi] / \langle 1 + \chi + \cdots + \chi^{N-1} \rangle
\]
is divisible by $1 + \chi$. Now we need the following lemma.
\begin{lem}
 Let $u \in 4\cdot \RhG^+$ be an element which maps to $0$ under the map $4 \cdot \RhG^+ \xra{\chi \mapsto -1} 8 \cdot \ZZ$. Then there exists $a \in 4 \cdot \RhG^-$ such that
 \[
  u = f \cdot a.
 \]
\end{lem}
\begin{proof}
 It is enough to prove the lemma for
 \[
  u = 4 \cdot (\chi^k + \chi^{-k}) + 8 \cdot (-1)^{k+1},
 \]
 where $k = 1,\ldots, N/2$. We have $u = (1+\chi) \cdot v$ where
 \begin{align*}
  v = \;& 4 \cdot (\chi^{k-1} - \chi^{k-2} + \cdots + (-1)^{k+1}) + \\ + \;& 4 \cdot (\chi^{-k} - \chi^{-k+1} + \cdots + (-1)^{k+1} \cdot \chi^{-1}).
 \end{align*}
 Setting $a = (1-\chi) \cdot v$ one easily verifies the desired equation as well as the fact that $a \in 4 \cdot \RhG^-$.
\end{proof}

Applying the above  lemma to $u = \wrho (x - \Sigma (y) - z \cdot
\sigma)$ we find an element $y' \in \wL^s_{4e+2} (G) \subset \sS^s
(L_N^{4e+1} (\alpha_k))$ such that $\wrho (y') = a$.

Using Lemma \ref{restriction-3-to-1} we see that $\eta \circ \Sigma
(y') = 0$. Further $\wrho (x - \Sigma (y+y') - z \cdot \sigma) = 0$.
Hence
\[
 x = \Sigma (y+y') + z \cdot \sigma \in \im (\Sigma) + \ZZ (\sigma).
    \]
To see that we obtain the direct sum it is enough to observe that
under the homomorphism
\[
 \sS^s (L_N^{4e+3} (\alpha_k)) \xra{\wrho} \QQ \RhG \xra{\chi \mapsto -1} \QQ
\]
the subgroup $\im (\Sigma)$ is mapped to $0$ (because $\wrho (\Sigma
(y)) = f \cdot \wrho (y)$) and that $\sigma$ is mapped to $8$ by
definition.

(3) This follows from the proof of the previous point.
\end{proof}

\subsection{Proof of Theorem \ref{thm2}}

\begin{lem} \label{lem_2-1}
Let $G = \ZZ_N$ with $N = 2^K \cdot M$, $M$ odd and $K \geq 1$. Then
\[
\bt_{4e-2} \circ \eta \circ \Sigma(x) = 0 \mbox{ for all } x \in
\sS^s(L_N^{4e-1} (\alpha_k)).
\]
\end{lem}
\begin{proof}
For $x \in \sS^s(L_N^{4e-1} (\alpha_k))$ we have
\[
\bt_{4e-2} \circ \eta \circ \Sigma(x) = \bt_{4e-2} \circ (p^N_2)^!
\circ \eta \circ \Sigma(x) = \bt_{4e-2} \circ \eta \circ \Sigma
\big((p^N_2)^!(x)\big).
\]
Therefore, it suffices to prove the statement for $N=2$. L{\'o}pez
de Medrano shows in \cite[chapter IV]{LdM(1971)} that $\sS^s
(L_2^{4e+1} (\alpha_k)) / \im(\Sigma) \cong \ZZ_2$ where the
isomorphism is induced by the map $\bt_{4e-2} \circ \eta$. In
particular, we have $\bt_{4e-2} \circ \eta \circ \Sigma (x) = 0$ for
all $x \in \sS^s(L_2^{4e-1} (\alpha_k))$.
\end{proof}

\begin{lem} \label{lem_2-2}
Let $G = \ZZ_N$ with $N = 2^K \cdot M$, $M$ odd and $K \geq 1$. Let
$\tau_N \in \wL^s_{4e}(\ZZ_N) \subset \sS^s (L_N^{4e-1} (\alpha_k))$
be defined by
\[
\Gsign (\tau_N) = 2^{\max\{ 4-K , 2 \}} \cdot (1 + \chi^2 + \cdots +
\chi^{N-2}) \in 4 \cdot \RhG^+.
\]
Then $\Sigma(\tau_N)$ has the following properties:
\begin{enumerate}
 \item $\wrho(\Sigma(\tau_N)) = 0$.
 \item $\bt_{2i} \circ \eta(\Sigma(\tau_N)) = 0$ if $i \neq 2e$.
 \item $\bt_{4e} \circ \eta(\Sigma(\tau_N)) \left\{ \begin{array}{ll} = 1 \in \ZZ_2 & K = 1 \\ \in \{ 2^{K-2} , 3 \cdot 2^{K-2} \} \subset \ZZ_{2^K} & K \geq 2 \end{array} \right..$
 \item $\bt_{(odd)} \circ \eta(\Sigma(\tau_N)) = 0$.
\end{enumerate}
\end{lem}
\begin{proof}
(1) $\wrho(\Sigma(\tau_N)) = f \cdot \Gsign(\tau_N) = 0$ because
\[
(1 + \chi) \cdot (1 + \chi^2 + \cdots + \chi^{N-2}) = 0 \in \QQ
\RhG.
\]

(2) $\res \circ \eta(\Sigma(\tau_N)) = \eta(\tau_N) = 0$ implies
$\bt_{2i} \circ \eta(\Sigma(\tau_N)) = 0$ if $i \notin \{ 2e-1 , 2e
\}$. Lemma \ref{lem_2-1} implies $\bt_{4e-2} \circ
\eta(\Sigma(\tau_N)) = 0$.

(3) First we consider the case $N = 2$. In \cite[chapter
IV]{LdM(1971)} L{\'o}pez de Medrano shows $\sS^s (L_2^{4e+1}
(\alpha_k)) / \im(\Sigma) \cong \ZZ_2$ where the isomorphism is
induced by the map $\bt_{4e-2} \circ \eta$. Let $x \in \sS^s
(L_2^{4e+1}(\alpha_k))$ be given by $\wrho(x) = 0$,
$\bt_{2i}(\eta(x)) = 0$ for $i \neq 2e$ and $\bt_{4e}(\eta(x)) = 1
\in \ZZ_2$. Because of the result of L{\'o}pez de Medrano there
exists an element $y \in \sS^s (L_2^{4e-1}(\alpha_k))$ with
$\Sigma(y) = x$. Since
\[
\eta(y) = \res\big(\eta\big(\Sigma(y)\big)\big) = \res(\eta(x)) = 0
\]
holds, we obtain $y \in \wL^s_{4e}(\ZZ_2) = \ZZ(\tau_2)$. Therefore,
we can write $y = m \cdot \tau_2$ with $m \in \ZZ$. We conclude
\[
m \cdot \bt_{4e} \circ \eta\big(\Sigma(\tau_2)\big) = \bt_{4e} \circ
\eta\big(\Sigma(y)\big) = \bt_{4e} \circ \eta(x) = 1 \in \ZZ_2
\]
and hence $\bt_{4e} \circ \eta(\Sigma(\tau_2)) = 1 \in \ZZ_2$.
	
Next we prove the statement for $N = 2^K$ by induction. We have
already checked the case $K = 1$. We proceed with the induction step
$K \to K+1$. We have to show $\bt_{4e} \circ
\eta(\Sigma(\tau_{2^{K+1}})) \in \{ 2^{K-1} , 3 \cdot 2^{K-1} \}
\subset \ZZ_{2^{K+1}}$ which is equivalent to
\[
\pr \circ \bt_{4e} \circ \eta(\Sigma(\tau_{2^{K+1}})) = 2^{K-1} \in
\ZZ_{2^K}
\]
with $\pr \colon \ZZ_{2^{K+1}} \to \ZZ_{2^K}, 1 \mapsto 1$. Notice
that
\[
(p^{2^{K+1}}_{2^K})^!(\tau_{2^{K+1}}) = \left\{ \begin{array}{ll}
\tau_2 & K = 1 \\ 2 \cdot \tau_{2^K} & K \geq 2 \end{array} \right.
\]
because
\[
\eta\big((p^{2^{K+1}}_{2^K})^!(\tau_{2^{K+1}})\big) =
(p^{2^{K+1}}_{2^K})^!\big(\eta(\tau_{2^{K+1}})\big) = 0 = \left\{
\begin{array}{ll} \eta(\tau_2) & K = 1 \\ \eta(2 \cdot \tau_{2^K}) &
K \geq 2 \end{array} \right.
\]
and
\begin{align*}
\wrho\big((p^{2^{K+1}}_{2^K})^!(\tau_{2^{K+1}})\big) & = 2^{\max\{ 3-K , 2 \}} \cdot (1 + \chi^2 + \chi^4 + \cdots + \chi^{2^{K+1}-2})\\
& = 2^{\max\{ 3-K , 2 \} + 1} \cdot (1 + \chi^2 + \chi^4 + \cdots + \chi^{2^K-2})\\
& = \left\{ \begin{array}{ll} \wrho(\tau_2) & K = 1 \\ \wrho(2 \cdot
\tau_{2^K}) & K \geq 2 \end{array} \right.
\end{align*}
hold. Using this result and the induction assumption we conclude
\begin{align*}
\pr \circ \bt_{4e} \circ \eta\big(\Sigma(\tau_{2^{K+1}})\big) & = \bt_{4e} \circ \eta\big(\Sigma\big((p^{2^{K+1}}_{2^K})^!(\tau_{2^{K+1}}) \big)\big)\\
& = \left\{ \begin{array}{ll} \bt_{4e} \circ \eta\big(\Sigma(\tau_2)\big) & K = 1 \\ \bt_{4e} \circ \eta\big(\Sigma(2 \cdot \tau_{2^K})\big) & K \geq 2 \end{array} \right.\\
& = 2^{K-1} \in \ZZ_{2^K}.
\end{align*}
It remains to show the statement for $N = 2^K \cdot M$ with $M > 1$.
We have the equation $(p^N_{2^K})^!(\tau_N) = M \cdot \tau_{2^K}$
because
\[
\eta\big((p^N_{2^K})^!(\tau_N)\big) = (p^N_{2^K})^!(\eta(\tau_N)) =
0 = \eta(M \cdot \tau_{2^K})
\]
and
\begin{align*}
\wrho\big((p^N_{2^K})^!(\tau_N)\big) & = 2^{\max\{ 4-K , 2 \}} \cdot (1 + \chi^2 + \chi^4 + \cdots + \chi^{N-2})\\
& = 2^{\max\{ 4-K , 2 \}} \cdot M \cdot (1 + \chi^2 + \chi^4 + \cdots + \chi^{2^K-2})\\
& = \wrho(M \cdot \tau_{2^K})
\end{align*}
hold. We conclude
\begin{align*}
\bt_{4e} \circ \eta\big(\Sigma(\tau_N)\big) & = \bt_{4e} \circ (p^N_{2^K})^! \circ \eta\big(\Sigma(\tau_N)\big)\\
& = \bt_{4e} \circ \eta\big(\Sigma\big((p^N_{2^K})^!(\tau_N)\big)\big)\\
& = M \cdot \bt_{4e} \circ \eta\big(\Sigma(\tau_{2^K})\big)\\
& \left\{ \begin{array}{ll} = 1 \in \ZZ_2 & K = 1 \\ \in \{ 2^{K-2}
, 3 \cdot 2^{K-2} \} \subset \ZZ_{2^K} & K \geq 2 \end{array}
\right..
\end{align*}

(4) We have $(p^N_M)^!(\tau_N) = 0$ because
\begin{align*}
\wrho\big((p^N_M)^!(\tau_N)\big) & = 2^{\max\{ 4-K , 2 \}} \cdot (1 + \chi^2 + \chi^4 + \cdots + \chi^{N-2})\\
& = 2^{\max\{ 4-K , 2 \} + K - 1} \cdot (1 + \chi + \chi^2 + \cdots + \chi^{M-1})\\
& = 0
\end{align*}
and $\eta\big((p^N_M)^!(\tau_N)\big) = (p^N_M)^!(\eta(\tau_N)) = 0$.
We conclude
\[
(p^N_M)^!\big(\eta\big(\Sigma(\tau_N)\big)\big) =
\eta\big(\Sigma\big((p^N_M)^!(\tau_N)\big)\big) = 0
\]
and hence $\bt_{(odd)} \circ \eta(\Sigma(\tau_N)) = 0$.
\end{proof}

\begin{lem} \label{lem_2-3}
Let $G = \ZZ_N$. If $t \in \ZZ$ satisfies $8 \cdot t \cdot f_k \in 4
\cdot \RhG$ then $N$ divides $4 \cdot t$.
\end{lem}
\begin{proof}
Notice that any element in $\RhG$ respectively in $\QQ\RhG$ can be
uniquely written as $\sum_{l=0}^{N-2} c_l \cdot \chi^{l \cdot k}$
with $c_l \in \ZZ$ respectively in. $\QQ$. We have
\[
f_k = 1 - \frac{2}{N} + \sum_{l=1}^{N-2} \Big( 2 - \frac{2+2 \cdot
l}{N} \Big) \cdot \chi^{l \cdot k}.
\]
Hence $8 \cdot t \cdot f_k \in 4 \cdot \RhG$ leads to an equation
\[
8 \cdot t \cdot \Big( 1 - \frac{2}{N} \Big) + \sum_{l=1}^{N-2} 8
\cdot t \cdot \Big( 2 - \frac{2+2 \cdot l}{N} \Big) \cdot \chi^{l
\cdot k} = \sum_{l=0}^{N-2} 4 \cdot c_l \cdot \chi^{l \cdot k}
\]
with $c_l \in \ZZ$. A comparison of the coefficients of $\chi^0$
gives $8 \cdot t \cdot ( 1 - 2/N ) = 4 \cdot c_0$ and hence $2 \cdot
t \cdot ( 1 - 2/N ) \in \ZZ$. This implies $4 \cdot t / N \in \ZZ$.
\end{proof}

\begin{proof}[Proof of Theorem \ref{thm2}]
(1) By Lemma \ref{lem_2-2} $\omega = 2^{\min\{ K , 2 \}} \cdot
\tau_N$ has the properties $\wrho(\Sigma(\omega)) = 0$ and
$\eta(\Sigma(\omega)) = 0$. This implies $\Sigma(\omega) = 0$ and
hence $\ZZ(\omega) \subseteq \ker(\Sigma)$. We want to prove
equality. Let $x \in \ker(\Sigma)$. We conclude $x \in
\wL^s_{4e}(\ZZ_N) \subset \sS^s (L_N^{4e-1} (\alpha_k))$ because
$\eta(x) = \res \circ \eta \circ \Sigma(x) = 0$. We have $f \cdot
\Gsign(x) = \wrho(\Sigma(x)) = 0$ and hence $\Gsign(x) \equiv 0 \;
\mod 1 + \chi^2 + \cdots + \chi^{N-2}$. This implies that $x$ is a
multiple of $\tau_N$. Since $\Sigma(x) = 0$, we conclude using Lemma
\ref{lem_2-2} that $x$ is a multiple of $\omega = 2^{\min\{ K , 2
\}} \cdot \tau_N$.

(2) Let $x_0 \in \sS^s(L_N^{4e+1}(\alpha_k))$. We want to show $x_0
\in \im(\Sigma) + \ZZ_2(\mu_{4e-2})$. We choose $y_1 \in \sS^s
(L_N^{4e-1}(\alpha_k))$ with $\eta(y_1) = \res(\eta(x_0))$. Then
$x_1 := x_0 - \Sigma(y_1)$ satisfies
\[
\res(\eta(x_1)) = \res(\eta(x_0)) - \res(\eta(\Sigma(y_1))) =
\eta(y_1) - \eta(y_1) = 0.
\]
If $M = 1$ then we set $y_2 := 0 \in \wL^s_{4e}(\ZZ_N) \subset
\sS^s(L_N^{4e-1}(\alpha_k))$. Otherwise we proceed as follows to
define $y_2$: Consider $\Sigma^{-1}((p^N_M)^!(x_1)) \in \sS^s
(L_M^{4e-1} (\alpha_k))$. We have $\Sigma^{-1}((p^N_M)^!(x_1)) \in
\wL^s_{4e}(\ZZ_M) \subset \sS^s(L_M^{4e-1}(\alpha_k))$ because
\[
\eta\big(\Sigma^{-1}\big((p^N_M)^!(x_1)\big)\big) =
\res\big(\eta\big((p^N_M)^!(x_1)\big)\big) =
(p^N_M)^!\big(\res\big(\eta(x_1)\big)\big) = 0.
\]
We choose $y_2 \in \wL^s_{4e}(\ZZ_N) \subset \sS^s (L_N^{4e-1}
(\alpha_k))$ with $(p^N_M)^!(y_2) = \Sigma^{-1}((p^N_M)^!(x_1))$. We
set $x_2 := x_1 - \Sigma(y_2) = x_0 - \Sigma(y_1 + y_2)$. Notice
that
\[
\res(\eta(x_2)) = \res(\eta(x_1)) - \res(\eta(\Sigma(y_2))) = 0 -
\eta(y_2) = 0
\]
and hence $\bt_{2i}(\eta(x_2)) = 0$ if $i \notin \{ 2e-1 , 2e \}$.
We have $\bt_{(odd)}(\eta(x_2)) = 0$ because
\[
(p^N_M)^!\big(\eta(x_2)\big) = (p^N_M)^!\big(\eta\big(x_1 -
\Sigma(y_2)\big)\big) = \eta\big((p^N_M)^!(x_1) -
\Sigma\big((p^N_M)^!(y_2)\big)\big) = \eta(0) = 0.
\]
Therefore, $\wrho(x_2)$ is of the shape $\wrho(x_2) = 8 \cdot t
\cdot f_k + r$ with $t \in \ZZ$ and $r \in 4 \cdot \RhGm$ (see
Proposition \ref{rho-formula-lens-sp}). Since $8 \cdot t \cdot f'_k
\in 4 \cdot \RhGp$ and
\begin{align*}
4 \cdot (\chi^l - \chi^{N-l}) = f \cdot \Big( & 4 \cdot (\chi^l + \chi^{N-l}) - 8 \cdot (\chi^{l+1} + \chi^{N-l-1}) + \cdots +\\
& (-1)^{N/2-l-1} \cdot 8 \cdot (\chi^{N/2-1} + \chi^{N/2+l}) +
(-1)^{N/2-l} \cdot 8 \cdot \chi^{N/2} \Big)
\end{align*}
for $1 \leq l \leq N/2-1$, there exists $z \in 4 \cdot \RhGp$ with
$f \cdot z = \wrho(x_2)$. We define $y_3 := \Gsign^{-1}(z) \in
\wL^s_{4e}(\ZZ_N) \subset \sS^s(L_N^{4e-1}(\alpha_k))$ and
\[
x_3 := x_2 - \Sigma(y_3) = x_0 - \Sigma(y_1 + y_2 + y_3).
\]
We have
\[
\wrho(x_3) = \wrho(x_2) - \wrho\big(\Sigma(y_3)\big) = \wrho(x_2) -
f \cdot z = 0.
\]
Since $\wrho\big((p^N_M)^!(x_3)\big) = 0$ holds and $\wrho \colon
\sS^s(L_M^{4e+1}(\alpha_k)) \to \QQ R_{\widehat \ZZ_M}^-$ is
injective, we obtain $(p^N_M)^!(x_3) = 0$. This implies
$(p^N_M)^!(\eta(x_3)) = 0$ and hence $\bt_{(odd)}(\eta(x_3)) = 0$.
Moreover, we have $\bt_{2i}(\eta(x_3)) = 0$ if $i \notin \{ 2e-1 ,
2e \}$ because
\[
\res(\eta(x_3)) = \res(\eta(x_2)) - \res(\eta(\Sigma(y_3))) = 0 -
\eta(y_3) = 0.
\]
Let $t_{4e-2} \in \{ 0 , 1 \}$ with $\bt_{4e-2}(\eta(x_3)) = t_{4e-2}$. We set
\[
x_4 := x_3 - t_{4e-2} \cdot \mu_{4e-2} = x_0 - \Sigma(y_1 + y_2 + y_3) - t_{4e-2} \cdot
\mu_{4e-2}
\]
and conclude $\wrho(x_4) = 0$, $\bt_{(odd)}(\eta(x_4)) = 0$ and
$\bt_{2i}(\eta(x_4)) = 0$ for $i \neq 2e$. By Proposition
\ref{rho-formula-lens-sp} we have $0 = \wrho(x_4) = 8 \cdot t \cdot
f_k + s$ where $s \in 4 \cdot \RhGm$ and $t \in \ZZ$ with $t \equiv
\bt_{4e}(\eta(x_4)) \; \mod 2^K$, $t \equiv 0 \; \mod M$. We know by
Lemma \ref{lem_2-3} that $2^{\max\{ K-2 , 0 \}}$ divides $t$. Hence
$\bt_{4e}(\eta(x_4))$ is a multiple of $2^{\max\{ K-2 , 0 \}}$. By
Lemma \ref{lem_2-2} there exists $n_\tau \in \ZZ$ with
$\bt_{4e}(\eta(x_4)) = n_\tau \cdot \bt_{4e}(\eta(\Sigma(\tau_N)))$.
Using again Lemma \ref{lem_2-2} we conclude that the element
\[
x_5 := x_4 - n_\tau \cdot \Sigma(\tau_N) = x_0 - \Sigma(y_1 + y_2 +
y_3 + n_\tau \cdot \tau_N) - t_{4e-2} \cdot \mu_{4e-2}
\]
satisfies $\wrho(x_5) = 0$, $\bt_{(odd)}(\eta(x_5)) = 0$ and
$\bt_{2i}(\eta(x_5)) = 0$ for all $i$. Hence $x_5 = 0$ and $x_0 =
\Sigma(y_1 + y_2 + y_3 + n_\tau \cdot \tau_N) + t_{4e-2} \cdot \mu_{4e-2} \in
\im(\Sigma) + \ZZ_2(\mu_{4e-2})$. This shows $\sS^s(L_N^{4e+1}(\alpha_k)) =
\im(\Sigma) + \ZZ_2(\mu_{4e-2})$. Notice that $\im(\Sigma) \cap \ZZ_2(\mu_{4e-2}) = 0$
because of Lemma \ref{lem_2-1} and $\bt_{4e-2}(\eta(a)) = 1$.
Therefore, we have $\sS^s(L_N^{4e+1}(\alpha_k)) = \im(\Sigma) \oplus
\ZZ_2(\mu_{4e-2})$.

(3) This follows from Lemma \ref{lem_2-1} and assertion (2).
\end{proof}

\begin{rem} \label{connection-to-BL}
The question when is a given element in the image of the suspension
map is called a {\it desuspension problem}. It was solved for $N =
2$ by L\'opez de Medrano in \cite{LdM(1971)} using the so-called
Browder-Livesay invariants, denoted $\sigma$. Depending on the
dimension one has $\sigma \co \sS (L^{4e-1}_2 (\alpha)) \ra \ZZ$, or
$\sigma \co \sS (L^{4e+1}_2 (\alpha)) \ra \ZZ_2$. In the dimension
$4e+1$ we have $\sigma = \bt_{4e-2} \circ \eta$. L\'opez de Medrano
also showed in \cite[section IV.4]{LdM(1971)} that in the dimensions
$4e-1$ the invariant $\sigma$ coincides with the $\wrho$-invariant
which is just an integer in the case $N = 2$.

Notice that the desuspension obstruction we obtain in the dimensions
$4e-1$ is given by sending $\chi \mapsto -1$ which corresponds to
passing to fake projective spaces via the transfer $(p_2^N)^!$
induced by the inclusion $\ZZ_2 \subset \ZZ_N$. Therefore our
obstructions to desuspension can be obtained by first applying
transfer and then taking the Browder-Livesay invariant.
\end{rem}

\section{Invariants of fake lens spaces} \label{sec:invariants}

Our main Theorem \ref{main-thm} gives an abstract calculation of the
structure sets of fake lens spaces. However, in order to obtain a
satisfactory classification one needs in addition to such a
calculation a good description of the invariants. In our case we
have a good description of the homomorphism into the free abelian
factor by the $\wrho$-invariant. But as already indicated in the
introduction the main theorem does not yield such a good description
for the torsion part.

We make a step towards a remedy of this deficiency. We still need to
make a choice of the splitting of the structure set into the free
abelian group and the torsion group as stated in Proposition
\ref{prop:sset-general}, see the proof of Theorem 5.1 in
\cite{Macko-Wegner(2008)}. However, once this choice is made, we obtain, in a special case, an isomorphism of the torsion part where the factors have certain geometric interpretation. It will be an
obstruction-theoretic description in a sense that we can interpret
certain factor if all the previous factors vanish.

The special case we treat is when $\alpha = \alpha_k$. Note that the
lens spaces $L^{2d-1}_N (\alpha_k)$ are obtained as joins $L^{1}_N
(\alpha_k) \ast L^{1}_N (\alpha_1) \ast \cdots \ast  L^{1}_N
(\alpha_1)$ and hence they contain sub-lens spaces $L^{2i-1}_N
(\alpha_k)$ for all $i \leq d$. This will be used in the description
of the invariants. Roughly speaking one examines the restrictions of
the simple homotopy equivalences which represent elements of $\sS^s
(L_N^{2d-1} (\alpha_k))_{\textup{tors}}$ to the inverse images of
these sub-lens spaces.

\begin{prop} \label{suspension-on-torsion}
Let $G = \ZZ_N$ with $N = 2^K \cdot M$, $M$ odd and $K \geq 1$, $e
\geq 2$. Then we have for the map
\[
 \Sigma \colon \sS^s (L_N^{4e-1} (\alpha_k)) \ra \sS^s (L_N^{4e+1} (\alpha_k))
\]
that:
\begin{enumerate}
 \item $2^{4 - \min\{K,2e\}} \cdot (1 + \chi^2 + \cdots + \chi^{N-2}) \in \im \; \big( \wrho \co \sS^s (L_N^{4e-1} (\alpha_k)) \ra \QQ \RhG \big)$.
 \item Let $\nu_{e} \in \sS^s (L_N^{4e-1} (\alpha_k))$ be such that
\[
\wrho (\nu_{e}) = 2^{4 - \min\{K,2e\}} \cdot (1 + \chi^2 + \cdots +
\chi^{N-2}).
\]
Then we have
\[
 \sS^s (L_N^{4e+1} (\alpha_k))_{\textup{tors}} \cong \im (\Sigma|_{\textup{tors}}) \oplus \ZZ_{2^{\min\{K,2e\}}} (\Sigma (\nu_{e})) \oplus \ZZ_2 (\mu_{4e-2})
\]
where $\mu_{4e-2} \in \sS^s (L_N^{4e+1} (\alpha_k))$ is defined by
$\wrho (\mu_{4e-2}) = 0$, $\bt_{(odd)} (\eta (\mu_{4e-2}) ) = 0$,
$\bt_{4e-2} (\eta (\mu_{4e-2}) ) = 1$ and $\bt_{2i} (\eta
(\mu_{4e-2}) ) = 0$ for $i \neq 2e-1$.
\end{enumerate}
\end{prop}

\begin{proof}
Let $l \in \ZZ$ be the smallest integer such that $2^l \cdot (1 +
\chi^2 + \cdots + \chi^{N-2})$ is in the image $\im ( \wrho \co
\sS^s (L_N^{4e-1} (\alpha_k)) \ra \QQ \RhG )$. Using the transfer
$(p_2^N)^{!}$ which sends $\chi \mapsto -1$ we conclude $2^l \cdot
2^{K-1} \cdot M \in 8 \cdot \ZZ$ and hence $l \geq 4 -K$. At the end
of the proof we show that $l = 4 - \min \{K,2e\}$. Before that we
prove the other statements.

Let $\nu_{e} \in \sS^s (L_N^{4e-1} (\alpha_k))$ be such that $\wrho
(\nu_{e}) = 2^l \cdot (1 + \chi^2 + \cdots + \chi^{N-2})$. Since
$\wrho (\Sigma (\nu_{e})) = f \cdot \wrho (\nu_{e}) = 0$ we obtain
$\Sigma (\nu_{e}) \in \sS^s (L_N^{4e+1} (\alpha_k))_{\textup{tors}}$.

We first show that
\[
\sS^s (L_N^{4e+1} (\alpha_k))_{\textup{tors}}    = \im
(\Sigma|_{\textup{tors}}) + \ZZ_{2^{\min\{K,2e\}}} (\Sigma
(\nu_{e})) + \ZZ_2 (\mu_{4e-2}).
\]
Abbreviate $\nu = \nu_{e}$ and $\mu = \mu_{4e-2}$ and let $y \in
\sS^s (L_N^{4e+1} (\alpha_k))_{\textup{tors}}$. Since we know $\sS^s
(L_N^{4e+1} (\alpha_k)) = \im (\Sigma) \oplus \ZZ_2 (\mu)$, we have
$y = \Sigma (x) + n_{\mu} \cdot \mu$ for some $x \in \sS^s
(L_N^{4e+1} (\alpha_k))$   .

Since $0 = \wrho (y) = f \cdot \wrho (x)$ we have that $\wrho (x) =
q \cdot (1 + \chi^2 + \cdots + \chi^{N-2})$ for some $q \in \QQ$.
Obviously $\bt_{(2)} (\eta (2^K \cdot x)) = 0$. If $M > 1$, then we
have $\bt_{(odd)} (\eta (2^K \cdot x)) = 0$ because
\begin{align*}
 \wrho ((p_M^N)^{!} (2^K \cdot x)) & = 2^K \cdot q \cdot (1 + \chi^2 + \cdots + \chi^{N-2}) = \\ & = 2^K \cdot q \cdot 2^{K-1} \cdot (1 + \chi + \chi^2 + \cdots + \chi^{M-1}) = 0
\end{align*}
and $(p_M^N)^{!} (2^K \cdot x) = 0$ since $\wrho$ yields an
injective homomorphism of $\sS^s (L_M^{4e-1} (\alpha_k))$. Hence
$\eta (2^K \cdot x) = 0$ and $2^K \cdot x \in \wL^s_{4e} (G)$ and
hence $2^K \cdot q \in \ZZ$.

We want to show $2^{-l} \cdot q \in \ZZ$. This is obviously true if
$q = 0$. If $q \neq 0$ then we can write $q = z \cdot 2^h$ with $h,z
\in \ZZ$ and $z$ odd. There exist $z',z'' \in \ZZ$ such that $(z'
\cdot z -1) \cdot 2^h = 2^l \cdot z''$. We obtain
\begin{align*}
 2^h \cdot (1 + \chi^2 + \cdots + \chi^{N-2}) & = z' \cdot \wrho (x) - z'' \cdot \wrho (\nu) = \\ & = \wrho (z' \cdot x - z'' \cdot \nu) \in \im \big( \wrho \co \sS^s (L_N^{4e-1} (\alpha_k)) \ra \QQ \RhG \big)
\end{align*}
which implies $h \geq l$ and $2^{-l} \cdot q = z \cdot 2^{h-l} \in
\ZZ$.

Set $n := 2^{-l} \cdot q \in \ZZ$. Notice $x - n \cdot \nu \in \sS^s
(L_N^{4e-1} (\alpha_k))_{\textup{tors}}$ because $\wrho (x - n \cdot
\nu) = \wrho (x) - n \cdot \wrho (\nu) = 0$.

We obtain
\begin{align*}
 y & = \Sigma (x) + n_{\mu} \cdot \mu = \\
 & = \Sigma (x - n \cdot \nu) + n \cdot \Sigma (\nu) + n_{\mu} \cdot \mu \in \im (\Sigma|_{\textup{tors}}) + \langle \Sigma (\nu) \rangle  + \ZZ_2 (\mu)
\end{align*}

A comparison of the orders of  the groups leads to $|\langle
\Sigma(x) \rangle | \geq 2^{\min\{K,2e\}}$. Since $\sS^s (L_N^{4e+1}
(\alpha_k))_{\textup{tors}}$ has no elements of order larger than
$2^{\min\{K,2e\}}$, we conclude $\langle \Sigma (x) \rangle \cong
\ZZ_{2^{\min\{K,2e\}}}$. Moreover, a comparison of the orders of the
groups shows
\[
 \sS^s (L_N^{4e+1} (\alpha_k))_{\textup{tors}}   = \im (\Sigma|_{\textup{tors}}) \oplus \ZZ_{2^{\min\{K,2e\}}} (\Sigma (\nu)) \oplus \ZZ_2 (\mu).
\]

It remains to show $l = 4 - \min\{K,2e\}$.

On one hand we have $2^{\min\{K,2e\}} \cdot \nu \in \ker (\Sigma)$
which implies $2^{\min\{K,2e\}} \cdot \nu = m \cdot \omega$ for some
$m \in \ZZ$. Therefore $2^{\min\{K,2e\}+l} - m \cdot 2^4 = 0$ and
hence $2^{\min\{K,2e\}+l-4} = m \in \ZZ$ and so $l \geq 4 -
\min\{K,2e\}$.

On the other hand let $\Delta := 2^{4-l} \cdot \nu - \omega$. This
has $\wrho (\Delta) = 0$, so we have $\Delta \in \sS^s (L_N^{4e-1}
(\alpha_k))|_{\textup{tors}}$. Since $\Sigma (\Delta) = 2^{4-l}
\cdot \Sigma (\nu)$ and $\im (\Sigma|_{\textup{tors}}) \cap \langle
\Sigma (\nu) \rangle = \{0\}$ we conclude $2^{4-l} \cdot \Sigma
(\nu) = 0$ and so $4-l \geq \min \{K,2e\}$ and $l \leq 4-\min
\{K,2e\}$.
\end{proof}

\begin{rem}
To avoid low-dimensional problems we need to make an ad-hoc
definition. Recall that
\[
 \sS^s (L^{5}_N (\alpha_k))_{\textup{tors}} \cong \ZZ_{2^{\max\{K,2\}}} \oplus \ZZ_2(\mu_{2}).
\]
Choose a generator of the first summand and denote it $\mu_4$.
\end{rem}

Recall from the previous section that we have an isomorphism on the
torsion parts $\Sigma \colon \sS^s (L_N^{4e+1}
(\alpha_k))_\textup{tors} \ra \sS^s
(L_N^{4e+3}(\alpha_k))_\textup{tors}$. Also note that $\nu_e \in
\sS^s (L^{4e-1} (\alpha_k))$ in Proposition
\ref{suspension-on-torsion} is not unique. Two choices of such
$\nu_e$ differ by a torsion element. Therefore a direct sum
decomposition of $\sS^s (L^{4e-1} (\alpha_k))$ into the free abelian
part and torsion part as in Proposition \ref{prop:sset-general}
determines $\nu_e$ uniquely by demanding that the projection on the
torsion part is $0$. So let us suppose that we have such a
decomposition for every $e \geq 2$. Then we obtain by induction:

\begin{cor} \label{cor:splitting-invariants}
Let $G = \ZZ_N$ with $N = 2^K \cdot M$, $M$ odd and $K \geq 1$, $d
\geq 3$. Then we have an isomorphism
\[
 (\bbr_{4i},\bbr_{4i-2})_i \co \sS^s (L_N^{2d-1} (\alpha_k))_{\textup{tors}} \cong \bigoplus_{i = 1}^c \ZZ_{2^{\min\{K,2i\}}} (\mu_{4i}) \oplus \bigoplus_{i = 1}^c \ZZ_2 (\mu_{4i-2})
\]
where $\mu_{4i} = \Sigma^{2(c-i)} (\nu_{i})$ when $i \geq 2$.
\end{cor}

\begin{rem}
The invariants $\bbr_{2i}$ are defined simply as the factors of the
isomorphism which follows by induction from Proposition
\ref{suspension-on-torsion}. The bar indicates that the invariants
$\bbr_{4i}$ differ from the invariants $\br_{4i}$ from our main
Theorem \ref{main-thm}. We prefer the new invariants since it
follows that they have a geometric interpretation as we
describe in more detail in the following remark.
\end{rem}

\begin{rem} \label{rem:splitting-invariants-4i}
The promised interpretation follows directly from Corollary
\ref{cor:splitting-invariants}. In detail, when $i \geq 2$, it
follows from the earlier discussion about splitting problems that
the invariants
\[
 \bbr_{4i} \co \sS^s (L_N^{4e+1} (\alpha_k)) \ra \ZZ_{2^{\min\{K,2i\}}} (\mu_{4i})
\]
can be described as follows. Let $h \co L(\alpha) \ra L_N
(\alpha_k)$ be a simple homotopy equivalence which represents $x \in
\sS^s (L_N^{4e+1} (\alpha_k))$. Suppose that $\bbr_{2j} (x) = 0$ for
$j > i$. Then we know that $h$ is homotopic to a simple homotopy
equivalence (still denoted $h$) which is split along $L_N^{4i-1}
(\alpha_k)$. Let $L_N^{4i-1} (\beta) := h^{-1} (L_N^{4i-1}
(\alpha_k))$. Then $h| \co L_N^{4i-1} (\beta) \ra L_N^{4i-1}
(\alpha_k)$ represents an element, say $y \in \sS^s (L_N^{4i-1}
(\alpha_k))$. Then $\bbr_{4i} (x)$ is the value of $y$ under the
composition
\[
 \sS^s (L_N^{4i-1} (\alpha_k)) \xra{\wrho} \QQ \RhG \xra {\chi \mapsto -1} \QQ \xra{:M:2^{3+\max\{0,K-2i\}}} \QQ.
\]
Notice that this value does not depend on the chosen decomposition
of $\sS^s (L^{4i-1} (\alpha_k))$ into the direct sum of the free
abelian part and the torsion part.
\end{rem}

\begin{rem} \label{rem:splitting-invariants-4i-2}
It follows from the results of the previous section that the
invariants $\bbr_{4i-2}$ also have a geometric description, when $i
\geq 2$. Namely, recall that in that case we have $\bbr_{4i-2} =
\br_{4i-2} = \bt_{4i-2}$ and after transferring via $(p_2^N)^!$
(which is an isomorphism on the $\ZZ_2 (\mu_{4i-2})$ summands) we
have $\bt_{4i-2} \circ \eta = \sigma$, where $\sigma$ is the
Browder-Livesay invariant.
\end{rem}

\begin{rem} \label{restriction-to-lens-spaces}
In the last two sections we have restricted ourselves to the simple
structure sets of $L^{2d-1}_N (\alpha)$ where $\alpha = \alpha_k$,
whereas in the main Theorem \ref{main-thm} this restriction  is not
present. To understand this, recall that we have first obtained the
calculation for $L^{2d-1}_N (\alpha_k)$ and then used the fact that
for any $\alpha$ there exists $k$ and a homotopy equivalence $h \co
L^{2d-1}_N (\alpha) \ra L^{2d-1}_N (\alpha_k)$. Such $h$ induces an
isomorphism on the simple structure sets. However, this induced
isomorphism has no good geometric description. It would certainly be
interesting to understand this problem.
\end{rem}

\small
\bibliography{lens-spaces}
\bibliographystyle{alpha}

\end{document}